\documentclass[10pt]{amsart}
\usepackage{fullpage}
\usepackage[french,english]{babel}
\usepackage{amsmath}
\usepackage{amsfonts}
\usepackage{amssymb}
\usepackage{amsthm}
\usepackage{graphicx}
\usepackage[latin1]{inputenc}
\usepackage[T1]{fontenc}
\usepackage{MnSymbol}
\usepackage[all]{xy}
\newtheorem{thm}{Theorem}[section]
\theoremstyle{definition}
\newtheorem{dfn}{Definition}[section]
\theoremstyle{remark}

\usepackage{color}
\newenvironment{dem}{\noindent{\textit{Proof.}}}{\begin{flushright}$\diamondsuit$\end{flushright}}
\theoremstyle{plain}
\newtheorem{lem}{\bf{Lemma}}[section]
\theoremstyle{plain}
\newtheorem{prop}{Proposition}[section]
\newtheorem{cor}{Corollary}[section]

\begin{document}
\author{}
\newcommand{\widering}[1]{\overset{\hbox{\smash{\lower1.333ex\hbox{$%
\displaystyle\mathring{}$}}}}{\wideparen{#1}}}
\def\Zz{{\mathbb{Z}}}
\def\Nn{{\mathbb{N}}}
\def\Cc{{\mathbb{C}}}
\def\Qq{{\mathbb{Q}}}
\def\End{{\mathrm{End}}}
\def\Hom{{\mathrm{Hom}}}
\def\Im{{\mathrm{Im}\,}}
\def\Aut{{\mathrm{Aut}}}
\def\Ker{{\mathrm{Ker}\,}}
\def\dim{{\mathrm{dim}\,}}
\def\Ext{{\mathrm{Ext}^1}}
\def\P1{{\mathbb{P}^1}}
\def\CalO{{\mathcal{O}}}
\def\CalF{{\mathcal{F}}}
\def\CalG{{\mathcal{G}}}
\def\CalE{{\mathcal{E}}}
\def\CalH{{\mathcal{H}}}
\def\CalV{{\mathcal{V}}}
\def\CalI{{\mathcal{I}}}
\def\CalH{{\mathcal{H}}}
\def\Coh{{\underline{\mathrm{Coh}}}}
\def\Lie{{\mathfrak{g}}}
\def\barp{{\underbar{p}}}
\def\deg{{\text{deg}}}
\def\tor{{\text{tor}}}
\def\canbun{{\vec{\omega}}}
\def\Vec{{\underline{\text{Bun}}}}
\def\Higgs{{\underline{\Lambda}_\mathbb{X}}}
\def\Higgsa{{\underline{\Lambda}_\mathbb{X}^\alpha}}
\def\Higgsb{{\underline{\Lambda}_\mathbb{X}^\beta}}
\def\Higgsab{{\underline{\Lambda}_\mathbb{X}^{\alpha+\beta}}}
\title[Loop crystal]{Higgs Bundles on weighted projective lines and loop crystals}
\author{Guillaume Pouchin}
\address{School of Mathematics, University of Edinburgh, JCMB King's Buildings, Edinburgh EH9 3JZ, United Kingdom.}
\email{g.pouchin@ed.ac.uk}
\thanks{The author gratefully acknowledges the financial support of EPSRC through the  grant EP/I02610x/1, as well as the University of Edinburgh. This work started during a stay at the Max Planck Institute, which is consequently thanked. The author would also like to thanks O. Schiffmann and I.Gordon for discussions and invaluable help, M. Chlouveraki for corrections on earlier versions, and the anonymous referee for corrections and suggestions.}
\maketitle
\begin{abstract}We consider the space of nilpotent Higgs bundles on a weighted projective line, as a global analog of the nilpotent cone. We show that it is pure, compute its dimension, and define geometric correspondences between irreducible components. We prove it by constructing a loop analog of a crystal, in the spirit of Kashiwara and Saito, for some corresponding loop Kac-Moody algebra.
\end{abstract}
\section{introduction}
The link between Kac-Moody algebras and geometry of quiver
representation is known since the work of Ringel (\cite{Ri}) and Lusztig (\cite{L2}), and has led to the theory of canonical bases of quantum algebras. More recently, Kapranov and Bauman-Kassel (\cite{BK}) have considered the Hall algebra of the category of coherent sheaves on the projective line and showed that some composition subalgebra is isomorphic to some positive part of the affine quantum algebra of $sl_2$. Schiffmann (\cite{Sc3}) generalized this situation to any weighted projective line as defined by Geigle and Lenzing in \cite{GL}, and showed that the Hall algebras obtained are affine versions of some quantum algebras. He used this setting to define canonical bases for some of these algebras. The category of coherent sheaves on weighted projective lines, or equivalently the category of bundles on $\P1$ with a parabolic structure, also appears in the study of the Deligne-Simpson problem, see \cite{CB1}.  

In the context of quivers, Lusztig considered a geometric
construction of (enveloping) Kac-moody algebras via constructible
functions on some nilpotent part $\Lambda_Q$ of the cotangent bundle of
the stack of representations of a quiver $Q$. He then defined
semicanonical bases using that geometry, whose elements are indexed by irreducible components of $\Lambda_Q$. Kashiwara and Saito (\cite{KS}) later considered natural correspondences between irreducible components to provide a geometric construction of the crystal associated to the semicanonical bases. 

The aim of this article is to develop Kashiwara and Saito's ideas in
the context of curves. We define the stack of nilpotent Higgs bundles $\Higgs$ on a weighted projective line $\mathbb{X}$, which was already consider in the case of the projective line in \cite{La2} (see also \cite{Gi}, \cite{Fa}, \cite{BD} for more general cases). We prove the following theorem, where $\Higgsa$ denote the substack parametrizing the Higgs bundles of class $\alpha\in K^+(\mathrm{Coh}(\mathbb{X}))$:
\begin{thm}
For any $\alpha$ in the Grothendieck group $K^+(\mathrm{Coh}_\mathbb{X})$, the stack $\underline{\Lambda}_\mathbb{X}^\alpha$ is pure of dimension $-\langle \alpha,\alpha \rangle$.
\end{thm}
Here the bracket $\langle.,.\rangle$ denotes the Euler form on the Grothendieck group.

We prove this by first constructing a stratification of this stack and then defining nice correspondences between strata which give rise to natural operators acting on the set irreducible components Irr($\Higgs$), indexed by indecomposable rigid coherent sheaves on the curve. These operators have the same properties as the classical crystal operators:
\begin{thm} For any $\alpha \in  K^+(\mathrm{Coh}_\mathbb{X})$ and any indecomposable rigid coherent sheaf $\CalI \in \mathrm{Coh}_\mathbb{X}$, there are operators:
\[
e_\CalI : \mathrm{Irr}(\underline{\Lambda}_\mathbb{X}^\alpha) \rightarrow \mathrm{Irr}(\underline{\Lambda}_\mathbb{X}^{\alpha+[\CalI]}),
\]
\[
f_\CalI : \mathrm{Irr}(\underline{\Lambda}_\mathbb{X}^\alpha) \rightarrow \mathrm{Irr}(\underline{\Lambda}_\mathbb{X}^{\alpha-[\CalI]}) \cup \{ 0\},
\]
and functions
\[
wt:\mathrm{Irr}(\underline{\Lambda}_\mathbb{X}) \rightarrow K^+(\Coh_\mathbb{X})
\]
and
\[
\epsilon_\CalI,\phi_\CalI: \underline{\Lambda}_\mathbb{X} \rightarrow \Zz
\]
with the following properties:
\begin{enumerate}
\item $wt(e_\CalI(Z))=wt(Z)+[\CalI]$ and $wt(f_\CalI(Z))=wt(Z)-[\CalI]$ if $f_\CalI(Z)\neq 0$,
\item $\epsilon_\CalI(e_\CalI(Z))=\epsilon_\CalI(Z)+1$ and $\phi_\CalI(e_\CalI(Z))=\phi_\CalI(Z)-1$,
\item $\epsilon_\CalI(f_\CalI(Z))=\epsilon_\CalI(Z)-1$  (resp. $\phi_\CalI(f_\CalI(Z))=\phi_\CalI(Z)+1$) if $\epsilon_\CalI(Z)\neq 0$, otherwise $f_\CalI(Z)=0$,
\item $\phi_\CalI(Z)=\epsilon_\CalI(Z)+<[\CalI],wt(Z)>$,
\item $Z'=f_\CalI(Z)$ if and only if $e_\CalI(Z')=Z$.
\end{enumerate}
Moreover for any two irreducible components $Z$ and $Z'$ in $\mathrm{Irr}(\underline{\Lambda}_\mathbb{X})$, there is a (finite) sequence of operators $u_j$ of the form $e_{\CalI_j}$ or $f_{\CalI_j}$, for $\CalI_j$ indecomposable rigid coherent sheaves, such that
\[
Z=Z_0 \xrightarrow{u_1} Z_1 \xrightarrow{u_2} Z_2 \xrightarrow{u_2} \cdots \xrightarrow{u_n} Z_n=Z'.
\]
\end{thm}
We define the \textit{loop crystal} to be this new combinatorial data, as it is an analog of the crystal structure in the case of quivers contructed in \cite{KS}. This is a first step towards a theory of crystals for loop Kac-Moody algebras, which, as in the usual Kac-Moody algebras case, is expected to carry very useful information about the representation theory of these algebras. The last property stated means that this loop crystal is connected. We will investigate the algebra of construcible functions on the space $\Higgs$ in a forthcoming paper. It is worth noting that our methods should apply to the quiver case, and give operators indexed by some other positive roots than just the simple roots. We also give some conjectures about the representations the loop crystal should be related to. Namely, in the case of affine Lie algebras, there is a description of integrable highest weight modules as semi-infinite limits of principal subspaces. Schiffmann constructed in \cite{Sc2} a canonical basis of these modules compatible with this description, which is conjecturally the same as the usual canonical basis. The loop crystal is expected to be the combinatorial structure carried by this basis.\\

The organization of the paper is the following: in the first part we introduce the definitions and notation of the space of coherent sheaves on a weighted projective line. In the second part we explain the algebraic stack structure on the space of Higgs bundles on these curves. The third part is devoted to the first properties of this space, giving a decomposition of irreducible components into their locally free part and their torsion part, then describing the irreducible components on the torsion part using results of Lusztig. The fourth part is the main part, where we introduce a stratification of the space, then define correspondences between irreducible components labelled by indecomposable rigid vector bundles. As a consequence, we prove that the space $\Higgs$ is Lagrangian. We also give some properties of the combinatorial data obtained, which we call \emph{loop crystal}. Finally, we prove in the last part that the loop crystal is connected as a colored graph. 
\section{Coherent sheaves on weighted projective lines}\label{projline}
In this section we recall the definitions of weighted projectives lines and their categories of coherent sheaves.
\subsection{Weighted projective lines}
In this section we recall the definition of weighted projective lines, as introduced in \cite{GL}, as well as the properties we will use. We keep the notations of \cite{Sc3}. Set $\barp=(p_1,p_2,\cdots,p_n)$ a tuple of $n$ positive integers. We will consider that $n\geq 3$, adding some $1$'s if needed.

We define the following:
\[
L(\barp)=\bigoplus_{i=1}^n \Zz \vec{x}_i /J,
\]
to be the quotient of the free abelian group generated by elements $\vec{x}_i$, $i=1\cdots n$, by the subgroup $J$ generated by the elements
\[
p_i\vec{x}_i-p_1\vec{x}_1
\]
for every $i\geq 2$.
Define $\vec{c}=p_1\vec{x}_1=p_i\vec{x}_i$. As $L(\barp)/\Zz \vec{c}\simeq \bigoplus_i \Zz/p_i\Zz$, $L(\barp)$ is an abelian group of rank $1$.

We also set $L^+(\barp)=\{l_1\vec{x}_1+ \cdots l_n\vec{x}_n\in L(\barp)\ | \ l_i\geq 0\}$ and define $\vec{x} \geq \vec{y}$ if $\vec{x}-\vec{y}\in L^+(\barp)$.\\
Define $S(\barp)=\Cc[X_1,\cdots, X_n]$ equipped with the $L(\barp)$-graduation given by $\mathrm{deg}(X_i)=\vec{x}_i$. Choose $n$ distinct points $\lambda_i$ in $\mathbb{P}^1(\Cc)$, with $\lambda_1=0$, $\lambda_2=\infty$ and $\lambda_3=1$. Consider the $L(\barp)$-homogeneous ideal $I_\lambda$ generated by polynomials $X_i^{p_i}-(X_1^{p_1}-\lambda_i X_2^{p_2})$ for $i \geq 3$.
We write $S(\barp,\lambda)=S(\barp)/I_\lambda$, still graded by $L(\barp)$.

Define the weighted projective line:
\[
\mathbb{X}=\mathbb{X}_{\barp,\lambda}:= \mathrm{Specgr}S(\barp,\lambda),
\]
together with its Zariski topology. The closed points in $\mathbb{X}$ are of two kinds:\\
-the ordinary points, corresponding to the ideal generated by $F(x_1^{p_1},x_2^{p_2})$, where $F$ is a prime homogeneous polynomial,\\
-the exceptional points $\lambda_1, \cdots,\lambda_n$, corresponding to the ideals generated by $x_i$.

We set $p=lcm(p_1,\cdots,p_n)$, and define the degree of a closed point by $\deg(x)=p$ if $x$ is ordinary and $\deg(\lambda_i)=\frac{p}{p_i}$. We define the \emph{genus} of $\mathbb{X}$ to be
\[
g_\mathbb{X}:=1+\frac{1}{2}\left( (n-2)p-\sum_{i=1}^n\frac{p}{p_i}\right).
\]
We associate to this data the star-shaped Dynkin diagram with $n$ branches of lengths $p_i-1$:
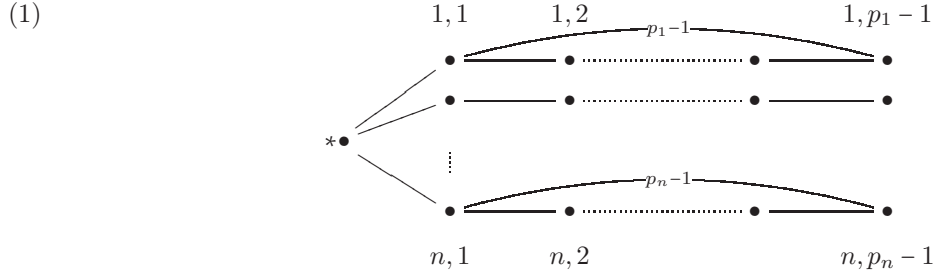
\begin{equation}\label{quiver}
\xymatrix@R=5pt{
 & 1,1 &1,2 & & & 1, p_1-1\\
  & \bullet \ar@{-}[r] \ar@/^1pc/ @{-}[rrrr]|{p_1-1}& \bullet \ar@{.}[rr] & & \bullet \ar@{-}[r] &\bullet \\
  & \bullet \ar@{-}[r] & \bullet \ar@{.}[rr]& & \bullet \ar@{-}[r] &\bullet \\
 \ast\bullet \ar@{-}[uur] \ar@{-}[ur]\ar@{-}[ddr]& \ar@{.}[d]& & & & \\
 &&&&&\\
& \bullet \ar@{-}[r] \ar@/^1pc/ @{-}[rrrr]|{p_n-1}& \bullet \ar@{.}[rr] & & \bullet \ar@{-}[r] &\bullet \\
 & n,1 &n,2 & & & n, p_n-1\\
}
\end{equation}
We denote by $J$ the set of vertices, and by $J_i$, $1\leq i \leq n$ the set of vertices corresponding to each branch of the diagram, so that $J=\bigcup_{i=1}^n J_i \cup \{*\}$.\\
As in the case of the projective line, we have a basis of the topology given by open subsets $D_f$, where $f$ is any homogeneous element of $S(\barp, \lambda)$, and:
\[
D_f=\{\mathfrak{p} \in \mathbb{X} \ | \ f\in \mathfrak{p}\}
\]
The structure sheaf $\CalO_\mathbb{X}$ of $\mathbb{X}$ is defined as the sheaf associated to the presheaf, $D_f \mapsto S(\barp,\lambda)_f$, where $S(\barp,\lambda)_f$ is the usual localization. We write $\CalO_\mathbb{X}$-Mod for the category of graded $\CalO_\mathbb{X}$-modules on $\mathbb{X}$. There is also a description via Serre's theorem:
\begin{thm}(Serre) (see \cite{GL}) The category of $\CalO_\mathbb{X}$ is equivalent to the category of $L(\barp,\lambda)$-graded modules over $S(\barp,\lambda)$ modulo finite length modules:
\[
\CalO_\mathbb{X} \simeq S(\barp,\lambda)-\mathrm{grmod}/S(\barp,\lambda)-\mathrm{grmod}_0,
\]
\end{thm}
where $S(\barp,\lambda)-\mathrm{grmod}_0$ are the $S(\barp,\lambda)$-modules of finite length.
Thanks to this theorem, we can shift objects in this category by elements of $L(\barp,\lambda)$. For any element $\vec{x} \in L(\barp,\lambda)$, we denote by $\CalO_\mathbb{X}(\vec{x})$ the $\CalO_\mathbb{X}$-module $\CalO_\mathbb{X}$ shifted by $\vec{x}$. Then for any $\CalO_\mathbb{X}$-module $\mathcal{M}$, we denote by $\mathcal{M}(\vec{x})$ the tensor product $\mathcal{M}\otimes_{\CalO_\mathbb{X}} \CalO_\mathbb{X}(\vec{x})$.
\subsection{Coherent sheaves}
The category $\mathrm{Coh}_\mathbb{X}$ of coherent sheaves on $\mathbb{X}$ is the category whose objects are $\mathcal{M} \in\CalO_\mathbb{X}$-Mod with the following property: there exists an open covering $\{U_i\}$ of $\mathbb{X}$ such that on each $U_i$ there is an exact sequence of $\CalO_{U_i}-$modules
\[
\bigoplus_{s=1}^N \CalO_\mathbb{X}(\vec{z}_s)_{|U_i} \rightarrow \bigoplus_{t=1}^M\CalO_\mathbb{X}(\vec{y}_t)_{|U_i} \rightarrow \mathcal{M}_{|U_i} \rightarrow 0
\]
This category shares many properties with the category of coherent sheaves on $\mathbb{P}^1$. We list the ones we will use later. We set $\canbun=(n-2)\vec{c} -\sum_{i=1}^n \vec{x}_i$.
\begin{prop}\cite{GL} We have the following properties:
\begin{enumerate}
\item The category $\Coh_\mathbb{X}$ is abelian and hereditary, i.e. for any $M,N\in\Coh_\mathbb{X}$, we have $\mathrm{Ext}^i(M,N)=0$ for $i>1$,
\item For any $M,N \in \Coh_\mathbb{X}$, the spaces $\Hom(M,N)$ and $\Ext(M,N)$ are finite dimensional,
\item (Serre duality) For any $\CalF,\CalG\in \Coh_\mathbb{X}$ there is a canonical isomorphism:
\[
\Ext(\CalF,\CalG)\simeq \Hom(\CalG,\CalF(\canbun))^*
\]
\end{enumerate}
\end{prop}

We denote by $\xi$ the point of $\mathbb{X}$ corresponding to the ideal $(0)$. The localisation $\CalO_{\mathbb{X},\xi}$ is Morita equivalent to the field $\Cc(\frac{X_2^{p_2}}{X_1^{p_1}})$. For a coherent sheaf $\CalF$, the fiber $\CalF_\xi=\CalF \otimes_{\CalO_\mathbb{X}}\CalO_{\mathbb{X},\xi}$ is a vector space over this field, and the dimension of this vector space is called the \textit{rank} of the sheaf $\CalF$, denoted by $\mathrm{rk}(\CalF)$.

Sheaves of rank $0$ are called \textit{torsion sheaves}. For such a sheaf $\CalF$ there are finitely many points $x\in\mathbb{X}$ where the fiber $\CalF_x$ is non zero: the set of such points is called the support of the sheaf $\CalF$, denoted by $\mathrm{supp}(\CalF)$. Any torsion sheaf splits in a direct sum of sheaves with disjoint supports.\\
The category of sheaves supported at an ordinary point $x$ is equivalent to the category of nilpotent representations of the Jordan (i.e. one loop) quiver. We denote by $\CalO_x^{(d)}$ the unique indecomposable sheaf of degree $d\delta$, and any sheaf supported on $x$ is isomorphic to a direct sum of sheaves of this type. \\
The category of sheaves supported at an exceptional point $\lambda_i$ has $p_i$ simples denoted by $S_j^{(i)}$, where $j\in \Zz/p_i\Zz$. The indecomposable sheaves are of the form $S_j^{(i)}(l)$, for $l$ a positive integer and $j\in \Zz/p_i\Zz$; $S_j^{(i)}(l)$ is uniquely determined by its head $S_j^{(i)}$ and its length $l$. 

A sheaf $\CalF$ such that there exists a covering $\{U_i\}$ with $\CalF_{|U_i}$ isomorphic to some $\bigoplus_{i=1}^N \CalO_\mathbb{X}(\vec{y}_i)_{|U_i}$ is called a \textit{vector bundle}.

Each coherent sheaf $\CalF$ has a (non-canonical) splitting $\CalF \simeq\CalF^\mathrm{vec} \oplus \CalF^\mathrm{tor}$, where $\CalF^\mathrm{vec}$ is a vector bundle and $\CalF^\mathrm{tor}$ is a torsion sheaf (see \cite{GL}).

We denote by $K(\mathbb{X})=K(\mathrm{Coh}_\mathbb{X})$ the Grothendieck group of $\mathrm{Coh}_\mathbb{X}$, and by $K^+(\mathrm{Coh}_\mathbb{X})$ its positive semigroup, whose elements $\alpha$ are such that there exists a coherent sheaf $\CalF$ with $[\CalF]=\alpha$ (where $[\CalF]$ denotes the class of $\CalF$ in $K(\mathrm{Coh}_\mathbb{X})$). The Grothendieck group $K(\mathrm{Coh}_\mathbb{X})$ is equipped with the Euler form: if $\alpha,\beta\in K^+(\mathrm{Coh}_\mathbb{X})$ then
\[
<\alpha,\beta>=\dim\Hom(\CalF,\CalG)-\dim\Ext(\CalF,\CalG), \ \mathrm{for\ any\ }\CalF,\CalG \ \mathrm{such\ that}\ [\CalF]=\alpha,\ [\CalG]=\beta.
\]
Remark that we also have for any $i$, $1\leq i\leq n$ and for any ordinary point $x\in \mathbb{X}$, we have in the Grothendieck group the equality
\[
\sum_{j=0}^{p_i-1} [S^{(i)}_{j}]=[\mathcal{O}_x].
\]
We also define the \emph{virtual genus} of the curve $\mathbb{X}$ to be 
\[
g_\mathbb{X}:= 1+ \frac{p}{2}(n-2-\sum_{i=1}^{n}\frac{1}{p_i}).
\]
It is easy to check that $g_\mathbb{X}<1$ if and only if the associated diagram is finite Dynkin, and $g_\mathbb{X}=1$ if and only if the associated diagram is affine Dynkin.
\subsection{Loop Kac-Moody algebras} In this section we follow the notations from (\cite{Sc3}).
The graph (\ref{quiver}) has an associated (symmetric, irreducible) Cartan matrix $C_\mathbb{X}=(a_{i,j})_{i,j\in J}$, where $a_{i,j}$ is such that $a_{i,i}=2$ and for $i\neq j$ it is equal to $-1$ if $i$ and $j$ are connected and $0$ otherwise. There is an associated (complex) Kac-Moody algebra $\mathfrak{g}_\mathbb{X}$ (see \cite{K}). We have
\[
\mathfrak{g}_\mathbb{X}=\mathfrak{n}_+ \oplus \mathfrak{h} \oplus \mathfrak{n}_-,
\]
and we denote by $e_i$, $i\in J$ (resp. $f_i$) the generators of the positive part $\mathfrak{n}_+$ (resp. negative part $\mathfrak{n}_-$). We denote by $\Delta$ the root system of $\mathfrak{g}_\mathbb{X}$, by $Q$ its root lattice, and by $\alpha_{i,j}$, $1\leq j < p_i$, $\alpha_*$ the simple roots corresponding to the vertices of the graph (\ref{quiver}). The lattice $Q$ is equipped with an invariant bilinear form $(.,.)$, such that $(\alpha_i,\alpha_j)=a_{i,j}$ for $i,j\in J$. \\
The \emph{loop} Lie algebra denoted by $\widehat{\mathfrak{g}}_\mathbb{X}$, is as a vector space:
\[
\widehat{\mathfrak{g}}_\mathbb{X}:= \mathfrak{g}_\mathbb{X} \otimes_\Cc \Cc[t,t^{-1}] \oplus \Cc \bf{c}.
\]
The Lie bracket between two elements $x\otimes t^l$ and $y\otimes t^k$, where $x$ and $y$ are in $\mathfrak{g}_\mathbb{X}$ is then given by
\[
[x\otimes t^l, y\otimes t^k]= [x,y]\otimes t^{l+k}+ l\delta_{l,-k}(x,y)\bf{c}.
\]
This Lie algebra is naturally equipped with a $\widehat{Q}=Q\oplus \Zz\delta$ grading, where the degree of $xt^l$ is the degree of $x$ plus $l\delta$, and the degree of $c$ is zero. The bilinear form on $Q$ is extended to $\widehat{Q}$ by defining $(\alpha,\delta)=0$ for any $\alpha\in \widehat{Q}$.\\
We also introduce the elements $\alpha_{i,0}\in \widehat{Q}$ for any $1\leq i \leq n$ defined by $\alpha_{i,0}=\delta-\sum_{j=1}^n\alpha_{i,j}$.\\
Denote by $\mathfrak{p}\subset \mathfrak{g}_\mathbb{X}$ the maximal parabolic subalgebra corresponding to the node $*$, and by $\mathfrak{l}$ and $\mathfrak{m}$ its Levi and nilpotent radical. Decompose the Levi component as $\mathfrak{l}=\mathfrak{n}_\mathfrak{l}^+ \oplus \mathfrak{h}\oplus \mathfrak{n}_\mathfrak{l}^-$ and define
\[
\widehat{\mathfrak{n}}=\mathfrak{n}_\mathfrak{l}^++t\mathfrak{l}[t] \oplus \mathfrak{m}[t,t^{-1}] \subseteq \widehat{\mathfrak{g}}.
\]
We consider an extension $\mathcal{L}\mathfrak{g}_\mathbb{X}$ of $\widehat{\mathfrak{g}}_\mathbb{X}$, given by the generators $e_{i,l}$, $f_{i,l}$ and $h_{i,l}$, for $i\in J$ and $l\in \Zz$, and relations:
\[
\begin{array}{cc}
[h_{i,l},h_{j,l'}]=l\delta_{l,-l'}a_{l,l'}\bf{c}, & [e_{i,l},f_{j,l'}]=\delta_{i,j}h_{i,l+l'} + l\delta_{l,-l'}\bf{c},\\

[h_{i,l},e_{j,l'}]=a_{i,j}e_{j,l+l'}, & [h_{i,l},f_{j,l'}]=-a_{i,j}f_{j,l+l'},\\

[e_{i,l+1},e_{j,l'}]=[e_{i,l},e_{j,l'+1}], & [f_{i,l+1},f_{j,l'}]=[f_{i,l},f_{j,l'+1}],\\

[e_{i,l_1},[e_{i,l_2},[\cdots [e_{i,l_n},e_{j,l'}]]\cdots ]=0 & \mathrm{if \ }n=1-a_{i,j},\\

[f_{i,l_1},[f_{i,l_2},[\cdots [f_{i,l_n},f_{j,l'}]]\cdots ]=0 & \mathrm{if \ }n=1-a_{i,j}.
\end{array}
\]
There is a natural surjective morphism $\mathcal{L}\mathfrak{g}_\mathbb{X} \rightarrow \widehat{\mathfrak{g}}_\mathbb{X}$, given by $e_{i,l} \mapsto e_i\otimes t^l$, $f_{i,l}\mapsto f_i \otimes t^l$ and $h_{i,l} \mapsto h_i\otimes t^l$.
When the virtual genus $g_\mathbb{X}$ is less than $1$, then the algebras $\widehat{\mathfrak{g}}_\mathbb{X}$ and $\mathcal{L}\mathfrak{g}_\mathbb{X}$ are isomorphic (see \cite{Sc3}).

The root system $\widehat{\Delta}$ of the algebra $\mathcal{L}\mathfrak{g}_\mathbb{X}$  and its root lattice $\widehat{Q}$ are then:
\[
\widehat{\Delta}=(\Delta +\Zz\delta) \cup \Zz^*\delta\ \mathrm{and} \ \widehat{Q}=Q+\Zz\delta.
\]
Introduce the positive part as follows: first remark that for any branch $J_i$, $1\leq i \leq n$, the subalgebra $\mathfrak{g}^i$ generated by $\{e_{j,l},f_{j,l},h_{j,l},\ j\in J_i\}$ is isomorphic to $\widehat{sl}_{p_i-1}$. Denote by $\theta_i$ the highest root of the subalgebra $sl_{p_i-1}$ and by $f_{-\theta_{i},l}$ the element of $\mathfrak{g}^i$ corresponding to $f_{-\theta_i}t^l$, where $f_{-\theta_i}$ is the lowest weight vector of $sl_{p_i-1}$. Then denote by $\widehat{\mathfrak{n}}^i$ the subalgebra of $\mathfrak{g}^i\subseteq \mathcal{L}\mathfrak{g}$ generated by $\{ e_{j,l},\ j\in J_i,\ l\in \Zz\} \cup \{ f_{-\theta_i,l},\ l\in \Zz\}$. Then $\mathcal{L}\mathfrak{n}$ is defined to be the subalgebra generated by
\[
\bigcup_{i=1}^n \mathfrak{n}^i \cup \{e_{*,l},\ l\in \Zz\}\cup \{h_{*,l},\ l\in \Zz\}.
\]
When the virtual genus $g_\mathbb{X}$ is less or equal to $1$, the algebras $\widehat{\mathfrak{n}}$ and $\mathcal{L}\mathfrak{n}$ are known to be isomorphic, see \cite{Sc3}.
The positive part of the root lattice $\widehat{Q}^+$, corresponding to the subalgebra $\mathcal{L}\mathfrak{n}$, is then formed by the elements
\[
k\alpha_*+c_\delta \delta+\sum_{i=1}^n\sum_{j=1}^{p_i-1} k_{i,j} \alpha_{(i,j)}
\]
where $k>0$ or $k=0$ and $c_\delta>0$ or $k=c_\delta=0$ and $k_{i,j}\geq 0$ for any $i,j$.\\
We have the following, see for instance \cite{Sc1} lemma 2.1 or \cite{Sc3} 5.1:
\begin{lem}
There is a canonical isomorphism $h: K(\mathbb{X})\simeq \widehat{Q}$, which maps $K^+(\mathbb{X})$ to the positive root lattice $\widehat{Q}^+$ and compatible with the bilinear forms. The explicit isomorphism is given by:
\[
\begin{array}{ccccc}
h: & K(\mathbb{X}) & \rightarrow  & \widehat{Q} &\\
 & [S_j^{(i)}] & \mapsto & \alpha_{i,j} &\\
 & [S_j^{(i)}(l)] & \mapsto & \sum_{k=0}^{l-1}\alpha_{i,j-k} &\\
 & [\CalO_x]  & \mapsto & \delta &  (x \mathrm{\ ordinary\ point)}\\
 & [\CalO_\mathbb{X}] & \mapsto & \alpha_* &
\end{array} 
\]
\end{lem}
Introduce the following set:
\[
\Pi= \{\alpha_*+\sum_{i=1}^n\sum_{k=1}^{j_i}\alpha_{i,k}\ | \ 0\leq j_i <p_i \}
\]
Associated to an element $t$ in $\Pi$, written $t=\alpha_*+\sum_{i=1}^n\sum_{k=1}^{j_i}\alpha_{i,k}$, we define the following line bundle:
\[
\mathcal{L}_t:=\CalO_\mathbb{X}(\sum_{i=1}^n j_i \vec{x}_i).
\]
\section{The stack of Higgs bundles on $\mathbb{X}$}
In this section we recall definitions and fix notation for the stack of Higgs bundles on a weighted projective line.
\subsection{The algebraic stack of coherent sheaves}
First we recall the stack structure on $\mathrm{Coh}_\mathbb{X}$. \\
We have a decomposition
\[
\Coh_\mathbb{X}=\bigsqcup_{\alpha \in K^+(\mathrm{Coh}_\mathbb{X})} \Coh_\mathbb{X}^\alpha,
\]
where $\Coh_\mathbb{X}^\alpha$ is the stack classifying isomorphism classes of coherent sheaves of class $\alpha$.
This stack has a local presentation as follows (see \cite{Gr},\cite{Le} or \cite{Sc1}).

Let $\alpha \in K^+(\mathrm{Coh}_\mathbb{X})$ and $\CalE \in \Coh_\mathbb{X}$.
Define the following functor from the category of schemes over $\Cc$ to the
category of sets: for $\Sigma$ a scheme over $\Cc$,
\[
\mathrm{\underline{Hilb}}_{\CalE,\alpha}(\Sigma)=\{ \phi_\Sigma: \CalE
\boxtimes \CalO_\Sigma \twoheadrightarrow
\CalF\ | \ \CalF \text{ is a coherent $\Sigma$-flat sheaf, } 
 \CalF_\sigma\text{ is of class $\alpha$ for all closed point $\sigma \in \Sigma$} \}/\sim
\]
where two such maps are equivalent if they have same kernel.
This functor is representable by a projective scheme
$\mathrm{Hilb}_{\CalE,\alpha}$ (see \cite{Gr} when $\mathbb{X}=\mathbb{P}^1$, and \cite{Sc1} in the general case of a weighted projective line).

Fix an integer $n\in \Zz$ and define for $t\in\Pi$:
\[ 
d_t(n,\alpha)=\langle[\mathcal{L}_t(n\vec{c})],\alpha\rangle \ \mathrm{\ and\ } \ 
\CalE_n^\alpha=\bigoplus_{t\in \Pi} \Cc^{d_t(n,\alpha)}\otimes \mathcal{L}_t(n\vec{c})\mathrm{\ if\ all\ }d_t(n,\alpha)\geq 0.
\]
A map $\phi_\Sigma$ induces for each closed point $\sigma \in \Sigma$ linear maps
$\phi_{t*,\sigma}:\Cc^{d_t(n,\alpha)} \rightarrow
\Hom(\mathcal{L}_t(n\vec{c}),\CalF_\sigma)$.\\
Let us consider the subfunctor defined by:
\[
\Sigma \mapsto \{ (\phi_\Sigma: \CalE_n^\alpha \boxtimes \CalO_\Sigma \twoheadrightarrow
\CalF ) \in \mathrm{\underline{Hilb}}_{\CalE_n^\alpha,\alpha}(\Sigma) \ | \ \forall
\sigma \in \Sigma,\  \forall t \in \Pi,\ \phi_{t*,\sigma} :\Cc^{d_t(n,\alpha)} \simeq
\Hom(\mathcal{L}_t(n\vec{c}),\CalF_\sigma) \}/\sim.
\]
This subfunctor is representable by a smooth open quasiprojective
subscheme $Q_n^\alpha$ of $\mathrm{Hilb}_{\CalE_n^\alpha,\alpha}$ (see \cite{Le}).\\
The group $G_n^\alpha= \prod_{t\in \Pi} \Aut(\mathcal{L}_t(n\vec{c})^{d_t(n,\alpha)})$ acts
naturally on $\mathrm{Hilb}_{\CalE_n^\alpha,\alpha}$ and $Q_n^\alpha$.\\
The quotient stacks $\Coh_\mathbb{X}^{\alpha,\geq n}:=[Q_n^\alpha/G_n^\alpha]$ for $n \in \Zz$ are open
substacks of $\Coh_\mathbb{X}^\alpha$ which form
an atlas, i.e.: 
\[
\Coh_\mathbb{X}^\alpha=\bigcup_{n \in \Zz}\Coh_\mathbb{X}^{\alpha,\geq
  n}=\bigcup_{n\in \Zz} [Q_n^\alpha/G_n^\alpha]
\]
We also introduce the stack $\Vec_\mathbb{X}=\bigcup_\alpha \Vec_\mathbb{X}^\alpha$ of vector bundles on
$\mathbb{X}$, which is an open substack of $\Coh_\mathbb{X}$. We have an atlas given by the open substacks
\[
\Vec_\mathbb{X}^{\alpha,\geq n}=[U_n^\alpha/G_n^\alpha]\ \mathrm{,\ where\ }\ U_n^\alpha=\{ (\phi:\CalE_n^\alpha \twoheadrightarrow \CalF)\in Q_n^\alpha\ | \ \CalF \mathrm{\ is\ locally\ free}\}/\sim.
\]
\subsection{The stack of Higgs bundles}
We want to describe the cotangent stack $T^*(\Coh_\P1^\alpha)$; an
atlas is obtained by symplectic reduction of the varieties
$T^*(Q_n^\alpha)$.\\
First recall that the tangent space at a point $\phi:\CalE
\twoheadrightarrow \CalF$ of $Q_n^\alpha$ is canonically isomorphic to
$\Hom(\Ker \phi,\CalF)$ (see \cite{Le}).\\
The group $G_n^\alpha$ acts on $T^*Q_n^\alpha$ in a Hamiltonian fashion. Let $\Lie_n^\alpha$ be the Lie algebra of $G_n^\alpha$. The corresponding moment map $\mu_n: T^*Q_n^\alpha \rightarrow
(\Lie_n^\alpha)^*$ is described as follows: 
 over a point $z=(\phi,f)\in T^*Q_n^\alpha$ with $\phi:
 \mathcal{E}_n^\alpha \rightarrow \CalF$ and $f\in \Hom(\Ker
 \phi,\CalF)^*$, we have 
\[
\begin{array}{cccc}
\mu_n^\alpha(z): &\Lie_n^\alpha&\rightarrow &\Cc\\
 & g  & \mapsto & <f,g|_{\Ker \phi}>,
\end{array}
\]
where $\Lie_n^\alpha$ is identified with $\bigoplus_{t\in\Pi}\Hom(\mathcal{L}_t(n\vec{c})^{d_t(n,\alpha)},\CalF)$ by means
of the isomorphisms $\phi_{t*}: \Cc^{d_t(n,\alpha)} \simeq
\Hom(\mathcal{L}_t(n\vec{c}),\CalF)$).\\
We want to describe the subvariety $(\mu_n^\alpha)^{-1}(0) \subseteq
T^*Q_n^\alpha$. To do this, fix some point $\phi: \mathcal{E}_n^\alpha
\twoheadrightarrow \CalF$ in $Q_n^\alpha$ and write the short exact sequence:
\[
0 \rightarrow \Ker(\phi) \rightarrow \mathcal{E}_n^\alpha
\rightarrow \CalF \rightarrow 0.
\]
Applying the functor $\Hom(\_,\CalF)$ yields
\[
0 \rightarrow \Hom(\CalF,\CalF)\rightarrow
\Hom(\mathcal{E}_n^\alpha,\CalF) \rightarrow \Hom(\Ker\phi, \CalF) 
\rightarrow \Ext(\CalF,\CalF)  \rightarrow
\Ext(\mathcal{E}_n^\alpha,\CalF) \rightarrow \cdots
\]
Since $\phi$ belongs to $Q_n^\alpha$, we have $\langle\mathcal{L}_t(n\vec{c}),\CalF \rangle=d_t(n,\alpha)=\dim\Hom(\mathcal{L}_t(n),\CalF)$,and so $\Ext(\mathcal{E}_n^\alpha,\CalF)=0$. By dualizing, we get :
\[
0 \rightarrow \Ext(\CalF,\CalF)^* \rightarrow \Hom(\Ker\phi,\CalF)^*
\xrightarrow{a} \Hom(\mathcal{E}_n^\alpha,\CalF)^* \rightarrow \cdots
\]
One can easily check that the map $a$ is precisely the moment map $\mu_n$. So if
$(\phi,f)$ is in $(\mu_n^\alpha)^{-1}(0)$ then $f$ defines a unique
element in $\Ext(\CalF,\CalF)^*$, which, abusing notation, we still denote by $f$.\\
Serre duality gives a canonical isomorphism: 
$\Ext(\CalF,\CalF)^* \simeq \Hom(\CalF,\CalF(\canbun))$, where we write
$\CalF(\canbun)$ for $\CalF \otimes \CalO(\canbun)$.\\
We finally have :
\[
(\mu_n^\alpha)^{-1}(0)=\{ (\phi : \mathcal{E}_n^\alpha \twoheadrightarrow
\CalF,f) \in T^*Q_n^\alpha\  \vert \ f \in \Hom(\CalF,\CalF(\canbun))\}/\sim
\]
By symplectic reduction, the cotangent bundle stack of the
quotient stack $Q_n^\alpha/G_n^\alpha$ is the quotient
$[(\mu_n^\alpha)^{-1}(0)/G_n^\alpha]$. This gives us an atlas of
$T^*(\Coh_\mathbb{X})$:
\[
T^*(\Coh_\mathbb{X}^\alpha)=\bigcup_{n\in \Zz}
[(\mu_n^\alpha)^{-1}(0)/G_n^\alpha].
\]
We then can write the objects of $T^*\Coh_\mathbb{X}^\alpha$ as follows:
\[
T^*\Coh_\mathbb{X}^\alpha=\{ (\CalF,f) \ | \ \CalF \in \Coh_\mathbb{X}^\alpha, \ f\in \Hom(\CalF,\CalF(\vec{\omega})) \}.
\]
\subsection{The global nilpotent cone}
Let us now introduce the nilpotent part of the cotangent bundle:
\[
S_n^\alpha:=(\mu_n^\alpha)^{-1}(0)^\mathrm{nilp} =\{
(\phi:\mathcal{E}_n^\alpha \twoheadrightarrow \CalF,f)\in \mu_n^{-1}(0)\ | \ f \mathrm{\
nilpotent}\}
\]
where we say that $f$ is \emph{nilpotent} if there exists $m$ such that
\[
f((m-1)\canbun)\circ\cdots \circ f(\canbun)\circ f =0
\]
as an element of $\Hom(\CalF,\CalF(m\canbun))$.\\
The quotient stacks $\underline{\Lambda}_\mathbb{X}^{\alpha,\geq
  n}=[S_n^\alpha/G_n^\alpha]$ are closed substacks of
  $T^*\Coh_\mathbb{X}^{\alpha,\geq n}$, and form a compatible family with
  respect to the inductive system
  $T^*\Coh_\mathbb{X}^{\alpha,\geq n}$. They give rise in the
  limit to a closed substack
\[\underline{\Lambda}_\mathbb{X}^\alpha=\varinjlim
  [S_n^\alpha/G_n^\alpha]=\bigcup_{n\in \Zz} [S_n^\alpha/G_n^\alpha]
  \subseteq T^*\Coh_\mathbb{X}^\alpha.
\]
We will then write objects of the global nilpotent cone as follows:
\[
\underline{\Lambda}_\mathbb{X}^\alpha=\{ (\CalF,f) \ | \ \CalF \in \Coh_\mathbb{X}^\alpha, \ f\in \Hom(\CalF,\CalF(\vec{\omega}))\ ,  \ f\mathrm{\ is\ nilpotent}\}.
\]
\section{First properties of irreducible components of $\underline{\Lambda}_\mathbb{X}$}\label{irredpart}
In this section we describe how the classification of irreducible of the global nilpotent cone splits into its locally free part and its torsion part, and then we describe the irreducible component of the torsion part.
\subsection{Vector bundles and torsion sheaves}
We first separate the study of irreducible components $\underline{\Lambda}_\mathbb{X}$ into their torsion and locally free parts. For an element $\CalF \in \Coh_\mathbb{X}$ we denote by $\CalF^\tor$ its torsion part. Let $\alpha,\beta \in K^+(\Coh_\mathbb{X})$, with $\mathrm{rk}(\beta)=0$, and define
\[
\underline{\Lambda}_\mathbb{X}^{\alpha,\beta}=\{ (\CalF,f) \in\Lambda_\mathbb{X}^\alpha \quad | \quad [\CalF^\tor]=\beta\}.
\]
We also define the stack parametrizing isomorphism classes of objects:
\[
\underline{L}^{\alpha,\beta}=\left\{ (\CalV, \tau,f_1,f_2,f_3)\ \left|\ \begin{array}{c} \CalV\in\Coh_\mathbb{X}^{\alpha-\beta,0}, \tau \in \Coh_\mathbb{X}^{\beta},\ f_1 \in\Hom(\CalV,\CalV(\canbun)) \mathrm{\ nilpotent}\\ f_2 \in\Hom(\tau,\tau(\canbun)) \mathrm{\ nilpotent}, \ f_3\in \Hom(\CalV,\tau(\canbun))\end{array} \right. \right\},
\]
an isomorphism between two objects $(\CalV,\tau,f_1,f_2,f_3)$ and $(\CalV',\tau',f_1',f_2',f_3')$ being a couple of isomorphisms $(\psi_1,\psi_2)$ where $\psi_1:\CalV \simeq \CalV'$ and $\psi_2: \tau \simeq \tau'$ such that the following diagrams commute:
\[
\xymatrix{
\mathcal{V} \ar[r]^{\psi_1} \ar[d]^{f_1} &\mathcal{V}' \ar[d]^{f_1'} &
\tau \ar[r]^{\psi_2} \ar[d]^{f_2}& \tau' \ar[d]^{f_2'} & \mathcal{V}
\ar[r]^{\psi_1} \ar[d]^{f_3} & \mathcal{V}'\ar[d]^{f_3'}\\
\mathcal{V}(\canbun) \ar[r]^{\psi_1(\canbun)} & \mathcal{V}'(\canbun)& \tau(\canbun) \ar[r]^{\psi_2(\canbun)} &
\tau'(\canbun) & \tau \ar[r]^{\psi_2} & \tau'
}
\]
We have a natural diagram:
\begin{equation}\label{Dec}
\xymatrix{
 & \underline{L}^{\alpha,\beta} \ar[dl]_{\pi_1} \ar[dr]^{\pi_2} & \\
\underline{\Lambda}_\mathbb{X}^{\alpha,\beta} & & \underline{\Lambda}_\mathbb{X}^{\alpha-\beta,0} \times \underline{\Lambda}_\mathbb{X}^{\beta}
}
\end{equation}
where $\pi_1$ is defined from the functor of groupoids:
\[
\begin{array}{cccc}
\pi_1: &(\mathcal{V},\tau,f_1,f_2,f_3)& \mapsto & \Bigg( \mathcal{V}
\oplus \tau,\begin{pmatrix}f_1&0\\f_3&f_2\end{pmatrix}\Bigg) \\
 & (\psi_1,\psi_2) & \mapsto & \begin{pmatrix}\psi_1 & 0 \\ 0 & \psi_2 \end{pmatrix}
\end{array}
\]
and $\pi_2$ is defined from the functor:
\[
\begin{array}{cccc}
\pi_2: & (\mathcal{V},\tau,f_1,f_2,f_3) & \mapsto &
((\mathcal{V},f_1),(\tau,f_2)) \\
 & (\psi_1,\psi_2)& \mapsto & (\psi_1,\psi_2)
\end{array}
\]
\begin{lem}\label{buntor}
The map $\pi_1$ is an affine fibration and $\pi_2$ is a vector bundle.
Both $\pi_1$ and $\pi_2$ are of relative dimension $\langle \alpha-\beta,\beta\rangle$ and with connected fibers. This
induces a bijection between irreducible components:
\[
Irr(\underline{\Lambda}_\mathbb{X}^{\alpha,\beta}) \leftrightarrow
Irr(\underline{\Lambda}_\mathbb{X}^{\alpha-\beta,0}) \times
Irr(\underline{\Lambda}_\mathbb{X}^{\beta}).
\]
Moreover this correspondence preserves dimensions, i.e. if we have $Z
\leftrightarrow Z_1 \times Z_2$ under this correspondence, then $\dim
Z=\dim Z_1 + \dim Z_2$.
\end{lem}
\begin{dem}
The result is obvious for $\pi_2$ since $\dim\Hom(\mathcal{V},\tau)=\langle \alpha-\beta,\beta \rangle$.\\
We also point out that the map $\pi_1$ is non representable.\\
We introduce the following natural stack parametrizing objects:
\[
\underline{\mathcal{S}}^{\alpha,\beta}= \{ (\CalF,f)| \CalF \in \Coh_\mathbb{X}^{\alpha,\beta},\ f\in
\Hom(\CalF,\CalF(\canbun)), f|_{\CalF^\mathrm{tor}}=0, f(\CalF)\subseteq
\CalF^\mathrm{tor}(\canbun) \}
\]
and a morphism $\psi$ between objects $(\CalF,f)$ and $(\CalF',f')$
is an isomorphism $\psi:\CalF \simeq \CalF'$ such that the diagram
\[
\xymatrix{
 \CalF \ar[r]^\psi \ar[d]^f & \CalF' \ar[d]^{f'} \\
\CalF(\canbun) \ar[r]^{\psi(\canbun)} & \CalF'(\canbun)
}
\]
is commutative.\\
We have a natural map $\pi: \underline{\mathcal{S}}^{\alpha,\beta}
\rightarrow \Coh_\mathbb{X}^{\alpha,\beta}$, which makes
$\underline{\mathcal{S}}^{\alpha,\beta}$ a vector bundle over $\Coh_\mathbb{X}^{\alpha,\beta}$.\\
Define its pullback over $\underline{\Lambda}_\mathbb{X}^{\alpha,\beta}$:
\[
\underline{\tilde{\mathcal{S}}}^{\alpha,\beta}=\underline{\mathcal{S}}^{\alpha,\beta}\times_{\Coh_\mathbb{X}^{\alpha,\beta}}\underline{\Lambda}_\mathbb{X}^{\alpha,\beta}
\]
This is a vector bundle over $\underline{\Lambda}_\mathbb{X}^{\alpha,\beta}$ of rank
$\langle \alpha-\beta,\beta\rangle$.\\
We have a natural action of $\underline{\tilde{\mathcal{S}}}^{\alpha,\beta}$
on $\underline{L}^{\alpha,\beta}$ over $\underline{\Lambda}_\mathbb{X}^{\alpha,\beta}$
defined as follows. Take a point $P=(\mathcal{V},\tau, f_1,f_2,f_3) \in
\underline{\Lambda}_\mathbb{X}^{\alpha,\beta}$. The fiber of
$\underline{\tilde{\mathcal{S}}}^{\alpha,\beta}$ over $\pi_1(P)$ is by construction
canonically identified with 
\[\{ g \in \End(\mathcal{V} \oplus \tau)|
g(\tau)=0,\ g(\mathcal{V} \oplus \tau) \subseteq \tau \}.
\]
We define the action of such elements $g$ as follows:
\[
g.P=(\mathcal{V},\tau,f_1,f_2,f_3-gf_1+f_2g);
\]
this corresponds to the action of $(Id +g)$ by conjugation on
$\begin{pmatrix} f_1 & 0 \\ f_3 & f_2 \end{pmatrix}$.\\
As we have $\Aut(\mathcal{V}\oplus \tau)= (\Aut(\mathcal{V}) \times
\Aut(\tau))\rtimes \Hom(\mathcal{V},\tau)$, we can identify
$\underline{\Lambda}_\mathbb{X}^{\alpha,\beta}$ as the quotient of
$\underline{L}^{\alpha,\beta}$ by the action of
$\underline{\tilde{S}}^{\alpha,\beta}$.
\end{dem}
\subsection{The cyclic quiver $\mathcal{C}_p$}
In order to study the irreducible components of the torsion part, we need to introduce the nilpotent variety for the space of representations of the cyclic quiver $\mathcal{C}_p$, where $p$ is a positive integer with $p\geq 2$. The cyclic quiver $\mathcal{C}_p$ is the quiver with vertices indexed by $i\in \Zz/p\Zz$, and arrows $\phi_i:i\rightarrow i+1$, for $i\in \Zz/p\Zz$. A \textit{nilpotent representation} of $\mathcal{C}_p$ is the following data:
\begin{itemize}
\item a set of finite dimensional $\Cc$-vector spaces $(V_i)_{i\in\Zz/p\Zz}$ indexed by the vertices,
\item for any $i \in \Zz/p\Zz$, a $\Cc$-linear map $\phi_i:V_i \rightarrow V_{i+1}$,
\end{itemize}
such that there exists a positive integer $k$ such that for any $i\in\Zz/p\Zz$, we have
\[
\phi_{i+k}\circ \phi_{i+k-1} \circ\cdots \circ {\phi_{i+1}}\circ \phi_i=0,
\]
as an element of $\Hom(V_i,V_{i+k+1})$.\\
The space of nilpotent representations of the cyclic quiver forms an abelian category. Denote by $\mathrm{RepNilp}(\mathcal{C}_p)$ the stack parametrizing isomorphism classes of nilpotent representations of $C_p$/ We will typically denote by $(V,\phi)$ a representation of $\mathcal{C}_p$. The Grothendieck groupof this category is isomorphic to $\Zz^p$. The class of a representation in the Grothendieck group is given by its dimension vector $(\dim V_i)_{i\in\Zz/p\Zz}$, so the semi-Grothendieck group $K^+(\mathrm{RepNil}(\mathcal{C}_p))$ is isomorphic to $\Nn^p$. The Grothendieck group is equipped with the Euler form.\\
We have the decomposition:
\[
\mathrm{RepNilp}(\mathcal{C}_p)=\bigsqcup_{\alpha\in \Nn^p}\mathrm{RepNilp}^\alpha(\mathcal{C}_p),
\]
where $\mathrm{RepNilp}^\alpha(\mathcal{C}_p)$ is the space of nilpotent representations of dimension vector $\alpha$.\\
For any $i\in \Zz/p\Zz$ and $l\in \Nn$, define the \textit{cyclic multisegment} $[i;l)$ to be the image of the projection to $\Zz/p\Zz$ of the segment $[i'-(l-1),i']$ for any $i'\in\Zz$, $i'\equiv i (mod\ p)$. A cyclic multisegment is a finite linear combination $\textbf{m}=\sum_{i,l}a_{i,l}[i;l)$ with $a_{i,l} \in \Nn$.
The objects of the category $\mathrm{RepNilp}(\mathcal{C}_p)$ are in natural correspondence with cyclic multisegments. Under this correspondence, the simple representation with dimension $1$ at the vertex $i$ and $0$ elsewhere is mapped to $[i;1)$, and the indecomposable representation of length $l$ and head $S_i$ is mapped to the cyclic segment $[i;l)$. We will also say that a cyclic multisegment $\textbf{m}=\sum_{i,l}a_{i,l}[i;l)$ is \textit{aperiodic} if for any $l\in \Nn$ there exists at least one $i$ such that $a_{i,l}$ is zero.\\ 
Introduce the doubled quiver $\overline{\mathcal{C}}_p$: it is obtained from $\mathcal{C}_p$ by adding an arrow from $i$ to $i-1$ for any $i\in \Zz/p\Zz$. A representation of the doubled quiver $\overline{\mathcal{C}}_p$ will be denoted by $(V,\phi,\bar{\phi})=(V_i,\phi_i,\bar{\phi}_i)_{i\in \Zz/p\Zz}$, where $(V,\phi)$ is a representation of $\mathcal{C}_p$ and $(\bar{\phi}_i)_{i\in \Zz/p\Zz}$ are maps
\[
\bar{\phi}_i: V_i \rightarrow V_{i-1}.
\]
Remark that there is no nilpotent condition on the maps $\bar{\phi}$.\\
Fix a dimension vector $\alpha$, and introduce the cotangent stack, which is obtained by symplectic reduction:
\[
T^*\mathrm{RepNilp}(\mathcal{C}_p)^\alpha:= \{(V,\phi,\bar{\phi})\in \mathrm{Rep}(\overline{\mathcal{C}}_p)^\alpha\ | \ [\phi,\bar{\phi}]=0 \},
\]
where $[\phi,\phi']=([\phi,\phi']_i)_{i\in\Zz/p\Zz}$, and $[\phi,\phi']_i=\phi'_{i+1}\phi_i-\phi_{i-1}\phi_i'\in \End(V_i)$.\\
The \textit{nilpotent variety}, introduced by Lusztig (\cite{L1},\cite{L5}), is defined as:
\[
\underline{\Lambda}_{\mathcal{C}_p}^\alpha:= \{ (V,\phi,\bar{\phi})\in T^*\mathrm{Rep}(\mathcal{C}_p)^\alpha,\ (V,\phi,\bar{\phi}) \mathrm{\ is\ nilpotent } \}.
\]
\begin{thm}(Lusztig, \cite{L2},\cite{L3},\cite{L5})\\
The stack $\underline{\Lambda}_{\mathcal{C}_p}^\alpha$ is pure of dimension $-\langle \alpha,\alpha  \rangle$.
\end{thm}
Remark also that the space $\underline{\Lambda}_{\mathcal{C}_p}$ has a rotation automorphism $\textbf{r}$: $\textbf{r} ((V_i,\phi_i,\bar{\phi}_i)_i)=(V_{i+1},\phi_{i+1},\bar{\phi}_{i+1})$.
\begin{thm}(\cite{LTV},\cite{KS})\\
The irreducible components of the space $\underline{\Lambda}_{\mathcal{C}_p}=\bigsqcup_{\alpha \in \Nn^p}\underline{\Lambda}_{\mathcal{C}_p}^\alpha$ are in natural bijection with the aperiodic cyclic multisegments.
\end{thm}
\subsection{Ordinary torsion sheaves}\label{Ordinary}
 Fix $\alpha\in K^+(\mathbb{X})$ such that $\mathrm{rank}(\alpha)=0$.\\
We first study the case of the projective line without any weight. Define, for a positive integer $d$, $\mathcal{P}(d)$ to be the set of \textit{partitions} of $d$, i.e. $\nu \in \mathcal{P}(d)$ is a set of non-negative integers $(\nu_1,\nu_2, \cdots)$ such that $\sum_i \nu_i=d$ and $\nu_{i+1}\leq \nu_i$. The \emph{length} $l(\nu)$ of a partition $\nu$ is the biggest integer $l$ such that $\nu_l\neq 0$.
For a partition $\nu$ of $d$ define
\[
\CalO_\P1^\nu:=\{
\CalO_{x_1}^{(\nu_1)}\oplus \cdots \oplus
\CalO_{x_d}^{(\nu_d)}\ | \ x_i \mathrm{\ distinct\ points}\},
\]
where $\CalO_x^{(d)}$ is the
indecomposable torsion sheaf supported on $x$ of degree $d$. This is a smooth strata of $\Coh_\P1^{d\delta}$. Let $T_{\CalO_{\P1,\nu}}^*\Coh_\P1^{d\delta}$ be the conormal bundle to this strata.
\begin{lem}[\cite{La1},theorem 3.3.13]\label{lemmetorsion}
We have the following decomposition into irreducible components:
\[
\underline{\Lambda}_\P1^{d\delta}=\bigcup_{\nu \in \mathcal{P} (d)}
\overline{T_{\CalO_\P1^\nu}^*\Coh_\P1^{d\delta}}
\]
Each irreducible component has dimension $0$.
\end{lem}
Now we go back to the case of an arbitrary weighted projective line $\mathbb{X}$ and we define for any positive integer $k$ the stack of torsion Higgs bundles of ordinary support:
\[
\underline{\Lambda}_{\mathbb{X},\mathrm{ord}}^{k\delta}:= \{ (\CalF,f) \in \underline{\Lambda}_\mathbb{X}^{d\delta} \ | \ \forall i=1,\cdots, n, \ x_i \notin \mathrm{Supp}(\CalF)\}
\]
Now remark that by definition the sheaf $\CalO_\mathbb{X}$ restricted to the open subset $\mathbb{X}-\{\lambda_i\}_i$ is obviously the same as the sheaf $\CalO_{\mathbb{P}^1}$ restricted to the open subset $\mathbb{P}^1-\{\lambda_i\}_i$. Then the stacks of torsion sheaves on $\mathbb{X}-\{\lambda_i\}_i$ and $\mathbb{P}^1-\{\lambda_i\}_i$ are isomorphic. The stack of torsion sheaves on $\mathbb{P}^1-\{\lambda_i\}_i$ is easily seen to be a dense open substack of the stack of torsion sheaves on $\P1$. We then deduce the following:
\begin{lem}
The stack $\underline{\Lambda}_{\mathbb{X},\mathrm{ord}}^{k\delta}$ is pure of dimension $0$. Its irreducible components are indexed by partitions of $k$.
\end{lem}
For a partition $\lambda$ of a positive integer $k$, we will denote by $Z_\lambda$ the corresponding irreducible component of $\underline{\Lambda}_{\mathbb{X},\mathrm{ord}}^{k\delta}$.\\

\subsection{Exceptional torsion sheaves}\label{Except}
For an element $\alpha^{(i)}$ such that $\alpha^{(i)}=\sum_{j=0}^{p_i-1}l_{(i,j)}\alpha_{i,j}$, define the stack of Higgs bundles supported on $\lambda_i$ of class $\alpha^{(i)}$:
\[
\underline{\Lambda}_{\mathbb{X},\lambda_i}^{\alpha^{(i)}}=\{(\CalF,f) \in \underline{\Lambda}_\mathbb{X}^{\alpha^{(i)}}\ | \ \mathrm{supp}(\CalF)=\lambda_i \}
\]
The category of coherent sheaves supported at $x_i$ is equivalent to the category of nilpotent representations of the cyclic quiver $\mathcal{C}_{p_i}$ (see \cite{Sc2}). This correspondence maps indecomposables sheaves $S_j^{(i)}(l)$ to the indecomposable element $S_j(l)$. This equivalence induces an isomorphism of stacks, and so we also deduce an isomorphism between the cotangent stacks:
\[
\begin{array}{cccc}
G: & \underline{\Lambda}_{\mathbb{X},\lambda_i}^{\alpha_l^i} & \rightarrow & \Lambda_{\mathcal{C}_{p_i}}^{(l_{i,j})}\\
 & (\CalF,f) & \mapsto & (V,\phi,\phi').
\end{array}
\]
It satisfies the following properties:
\begin{itemize}
\item the map $\CalF \rightarrow (V,\phi)$ is deduced from the equivalence $\Coh_{\mathbb{X},x_i}^{\alpha_l^i} \simeq \mathrm{RepNil}_{\mathcal{C}_{p_i}}^{\alpha_l^i}$.
\item tensoring by $\canbun$ on $\Coh_{\mathbb{X},x_i}^{\alpha_l^i}$ corresponds via $\Phi$ to the diagram automorphism $\textbf{r}$ acting on $\mathrm{RepNil}_{\mathcal{C}_{p_i}}^{\alpha_l^i}$, i.e. $\textbf{r}((V_i,\phi_k)_{k\in \Zz/ p_i\Zz})=(V_{k+1},\phi_{k+1})_{k\in \Zz/ p_i\Zz}$.
\item the element $f\in \Hom(\CalF,\CalF(\canbun))$ is mapped via $\Phi$ to an element $\phi' \in \Hom(\phi,\textbf{r}_*^{-1}\phi)$, so that if $\phi=(V_k,\phi_k:V_k \rightarrow V_{k+1})_{k\in \Zz/ p_i\Zz}$, then the element $\phi'$ is some element:
\[
\phi'=(\phi_k':V_k \rightarrow V_{k-1})
\]
which is nilpotent and commutes with $\phi$, i.e. $\phi_{k-1}\phi_k'=\phi_{k+1}'\phi_k$ for $k\in \Zz /p_i\Zz$. So $(\phi,\phi')$ is actually an element of:
\[
\Lambda_{\mathcal{C}_{p_i}}^{(l_{i,j})}=\{(\phi,\phi') \in \mathrm{RepNil}(\overline{\mathcal{C}}_{p_i}^{\alpha_l^i})\ |\ (\phi,\phi') \mathrm{\ nilpotent}, [\phi,\phi']=0\}.
\]
\end{itemize}
Reciprocally, a nilpotent element $\phi'\in \Hom(\phi,\textbf{r}_*\phi)$ is easily seen to correspond to an element $f\in\Hom_\Cc(\CalF,\CalF(\canbun))$, and the commutativity implies that it is in fact a morphism of $\CalO_\mathbb{X}$-modules, which is nilpotent.\\
As a consequence, we have:
\begin{lem}
The stack $\underline{\Lambda}_{\mathbb{X},\lambda_i}^{\alpha_l^i}$ is isomorphic to $\Lambda_{\mathcal{C}_{p_i}}^{(l_{i,j})}$. Each irreducible component is of dimension $-\langle \alpha_l^i,\alpha_l^i\rangle$. The irreducible components are indexed by aperiodic cyclic $p_i$-multisegments of dimension $\alpha_l^i$.
\end{lem}
\subsection{Irreducible components of the stack of torsion Higgs bundles}\label{IrredTor}
In this subsection we prove that the space $\underline{\Lambda}_\mathbb{X}^\alpha$, with $\mathrm{rk}(\alpha)=0$, is pure of dimension $-\langle \alpha,\alpha \rangle$, and give a description of the irreducible components.
Define for an element $\alpha \in K^+(\Coh_\mathbb{X})$ such that $\mathrm{rk}(\alpha)=0$ the set
\[ 
W(\alpha)=\{(l,(l_{i,j})_{i=1\cdots n, j\in \Zz/ p_i\Zz})\ |\ l, l_{i,j}\in \Nn,\ l\delta +\sum_{i,j}l_{i,j} \alpha_{i,j}=\alpha \}.
\]
Any $\CalF \in \Coh_\mathbb{X}^\alpha$ is canonically decomposed as
\[
\CalF=\CalF_\delta \oplus \bigoplus_{i=1}^n \CalF_i,
\]
where $\mathrm{supp}(\CalF_\delta)\subseteq \mathbb{X}-\{\lambda_i\}_{i=1\cdots n}$ and for any $i=1,\cdots,n, \mathrm{supp}(\CalF_i)=\lambda_i$.

We introduce the following stratification of $\Higgsa$: for $w=(l,l_{(i,j)}) \in W(\alpha)$, define $\alpha_w^{(i)}:=\sum_jl_{(i,j)}\alpha_{i,j}$ and 
\[
\underline{\Lambda}_\mathbb{X}^{\alpha,w}=\{ (\CalF,f) \in \Higgsa\ |\  [\CalF_\delta]=l\delta, \ \forall i \ [\CalF_i]=\alpha_w^{(i)} \}
\]
For such sheaves $\CalF$, the elements $f$ split as a sum $f=f_\delta \oplus \bigoplus_i f_i$, where $f_\delta \in\Hom(\CalF_\delta,\CalF_\delta(\canbun))$ and $f_i\in\Hom(\CalF_i,\CalF_i(\canbun))$, all of which are nilpotent.

By definition we have
\[
\underline{\Lambda}_\mathbb{X}^{\alpha,w} \simeq \underline{\Lambda}_{\mathbb{X},\mathrm{ord}}^{l\delta} \times \prod_{i=1}^n\underline{\Lambda}_{\mathbb{X},\lambda_i}^{\alpha_w^{(i)}},\quad \mathrm{\ and\ }\quad \underline{\Lambda}_\mathbb{X}^\alpha=\bigsqcup_{w\in W(\alpha)} \underline{\Lambda}_\mathbb{X}^{\alpha,w}.
\]
Combining the preceding lemmas, we obtain the following description of the irreducible component of the torsion part. 
\begin{thm} Let $\alpha\in K^+(\Coh_\mathbb{X})$ such that $\mathrm{rk}(\alpha)=0$. We have a natural bijection:
\[
\mathrm{Irr}(\underline{\Lambda}_\mathbb{X}^\alpha) \leftrightarrow \bigsqcup_{(l,l_{(i,j)})\in W(\alpha)} \mathcal{P}(l) \times \mathrm{Irr}(\Lambda_{\mathcal{C}_{p_i}}^{\alpha_l^i})
\]
Each irreducible component is of dimension $-\langle \alpha,\alpha\rangle$.
\end{thm}

\section{The loop crystal}\label{Loop}
In this section we define natural correspondences between irreducible components, indexed by indecomposable rigid coherent sheaves.
\subsection{Geometric correspondences}
Fix an indecomposable rigid element $\mathcal{I}$ in the category $\mathrm{Coh}_\mathbb{X}$.
For a coherent sheaf $\CalF$ on $\mathbb{X}$, define 
\[
\mathrm{\epsilon}_\CalI(\CalF)=\mathrm{max}\{i\in \Nn\ | \ \mathrm{Inj}(\CalI^{\oplus i},\CalF)
\neq \emptyset\}
\]
where $\mathrm{Inj}(\CalI^{\oplus i},\CalF)$ is the set of injections in
$\Hom(\CalI^{\oplus i},\CalF)$.\\
Fix $\CalI$ an indecomposable rigid element of $\mathrm{Coh}_\mathbb{X}$, two non-negative integers $s$ and $n$, an element $\alpha\in K^+(\mathrm{Coh}_\mathbb{X})$ and define the following locally closed substacks:
\[
\underline{\Lambda}_{\CalI,s}^\alpha:= \{(\CalF,f) \in
\underline{\Lambda}^\alpha_\mathbb{X}\ |\  \mathrm{\epsilon}_\CalI(\Ker(f))=s\},
\]
\[
\underline{\Lambda}^\alpha_{\CalI,n,s}=\{(\CalF,f) \in
\underline{\Lambda}^\alpha_\mathbb{X}\ |\  \mathrm{\epsilon}_\CalI(\Ker(f))=s,\
\dim(\Hom(\CalI,\Ker(f)))=n \}.
\]
We have that:
\[
\underline{\Lambda}_{\mathbb{X}}^\alpha=\bigcup_{s\geq 0}\underline{\Lambda}_{\CalI,s}^\alpha,
\]
a finite union of constructible substacks, and
\[
\underline{\Lambda}_{\mathbb{X}}^\alpha=\bigcup_{s,n\geq 0}\underline{\Lambda}_{\CalI,n,s}^\alpha,
\]
a locally finite union of constructible substacks (i.e. the number of non-empty intersections of these strata with a substack of $\Lambda_\mathbb{X}$ of finite type is finite).\\ 
In the following write $\gamma=s[\CalI]$ and $\beta=\alpha-\gamma$, three elements of $K(\mathrm{Coh}_\mathbb{X})$.\\
Define the following stack:
\[
\underline{\mathcal{E}}^\alpha_{\CalI,n,s}=\{ (\CalF,f,i)\ | \ (\CalF,f) \in
\underline{\Lambda}^{\alpha}_{\CalI,n,s}, \  i \in
\mathrm{Inj}(\CalI^{\oplus s},\Ker(f)) \}
\]
representing the functor from the category of affine schemes to the
category of groupoids:
\[
\Sigma\mapsto \left\{ (\CalF,f,i) \ \left| \begin{array}{c} \CalF
  \text{ is a coherent $\Sigma$-flat sheaf on $\mathbb{X}\times \Sigma$, } f \in \Hom(\CalF, \CalF \otimes \CalO_{\Sigma\times\mathbb{X}}(\canbun)),\
  i:\mathcal{I}\boxtimes \CalO_{\Sigma\times \mathbb{X}}^s \rightarrow \Ker f,  \\
 \CalF_\sigma \text{ is of class $\alpha$ and $\mathrm{\epsilon}_\CalI(\Ker f_\sigma)=s$},\ i_\sigma: \CalI^s\hookrightarrow  \Ker f_\sigma  \text{ for all closed point $\sigma \in \Sigma$}\end{array} \right. \right\}
\]
where a morphism $\psi$ between objects $(\CalF,f,i)$ and
$(\CalF',f',i')$ is an isomorphism $\psi:\CalF\simeq \CalF'$ such that
the following diagrams commute:
\[
\xymatrix{
 \CalF \ar[r]^\psi \ar[d]^f &  \CalF'\ar[d]^{f'} & \CalI^s\boxtimes \CalO_\Sigma
 \ar[r]^i \ar[rd]^{i'} & \Ker f \ar[d]^\psi \\
\CalF(\canbun) \ar[r]^{\psi(\canbun)} & \CalF'(\canbun) &  & \Ker f'
}
\]
Now consider the following diagram:
\begin{equation}\label{diag}
\xymatrix{
 & \underline{\mathcal{E}}^{\alpha}_{\CalI,n,s}  \ar[ld]_{p_1} \ar[rd]^{p_2}& \\
\underline{\Lambda}^{\alpha}_{\CalI,n,s} & &
 \underline{\Lambda}^{\beta}
}
\end{equation}
where the maps $p_1$ and $p_2$ are defined on objects as follows:
\[
\begin{array}{cccc}
p_1: & (\CalF,f,i) & \mapsto & (\CalF,f) \\
p_2: & (\CalF,f,i)& \mapsto & (\CalF/i(\CalI^s),f|_{\CalF/i(\CalI^s)}),
\end{array}
\]
We begin with the following lemma:
\begin{lem}
$p_2(\underline{\mathcal{E}}^{\alpha}_{\CalI,n,s})\subseteq \underline{\Lambda}_{\CalI,n-s,0}^{\beta}$.
\end{lem}
\begin{dem}
Take an element $(\CalF,i,f)$ in
$\underline{\mathcal{E}}^{\alpha}_{\CalI,n,s}$ and define
$(\CalG,g)=p_2(\CalF,i,f)$. We have a short exact sequence:
\begin{equation}
\xymatrix{
0 \ar[r] & \CalI^s \ar[r]^i & \CalF \ar[r]^j & \CalG \ar[r] & 0
}
\end{equation}
As $\Ext(\CalI,\CalI)=0$, we also have
\begin{equation}
\xymatrix{
0 \ar[r] & \Hom(\CalI,\CalI^s) \ar[r]^{i' }& \Hom(\CalI,\CalF) \ar[r]^{j'} & \Hom(\CalI,\CalG) \ar[r] & 0
}
\end{equation}
Now the following diagram commutes for any $a\in\Hom(\CalI,\Ker f
)$
\begin{equation}
\xymatrix{
\CalI \ar[r] \ar[rd]^a & \CalG \ar[r]^g & \CalG(\canbun) \\
 & \CalF \ar[u]^j \ar[r]^f & \CalF(\canbun) \ar[u]^j,
}
\end{equation}
so we see that $f\circ a=0$ implies $g\circ j'(a)=0$, whence
$j'(\Hom(\CalI,\Ker f ))\subseteq \Hom(\CalI,\Ker g)$. In the
other way, if $g\circ j'(a)=0$, then as the kernel of the morphism
$\CalF(\canbun) \rightarrow \CalG(\canbun)$ is $\CalI(\canbun)^s$, we have that
$\Im(f\circ a)\subseteq \CalI(\canbun)^s$. The morphism $f\circ a$ lies
inside $\Hom(\CalI,\CalI(\canbun)^s)=(\Ext(\CalI,\CalI)^s)^*=0$, so that $a\in
\Hom(\CalI,\Ker f)$. We just proved that we have a short exact sequence:
\[
0 \rightarrow \Hom(\CalI,\CalI^s) \rightarrow \Hom(\CalI,\Ker
f) \rightarrow \Hom(\CalI,\Ker g) \rightarrow 0
\]
It follows that $\dim \Hom (\CalI,\Ker g)=n-s$.\\
Now we prove that there are no injections from $\CalI$ into
$\Ker g$. Assuming that such an injection $h$ exists, consider an element
$h'\in \Hom(\CalI,\Ker f)$ such that $j'(h')=h$. Define
$h''=h'\oplus i \in \Hom(\CalI\oplus \CalI^s,\Ker f)$. From the
commutative diagram: 
\begin{equation*}
\xymatrix{
0 \ar[r] & \CalI^s \ar[r]  & \CalI^{s+1}
\ar[r]^{pr}\ar[d]^{h''}&  \CalI\ar[d]^h \ar[r]& 0\\
 & & \CalF \ar[r]^j & \CalG & ,
}
\end{equation*}
we deduce that $pr(\Ker h'')=0$, i.e. $\Ker h''\subseteq \CalI^s$. But the restriction of
$h''$ to $\CalI^s$ is $i$, which is injective. We have thus proved that
$h''$ is an injection from $\CalI^{s+1}$ into $\Ker f$, which is
impossible since $\mathrm{\epsilon}_\CalI(\Ker f)=s$.
\end{dem}
The diagram (\ref{diag}) is then refined to the diagram:
\begin{equation}\label{diag2}
\xymatrix{
 & \underline{\mathcal{E}}^{\alpha}_{\CalI,n,s}  \ar[ld]_{p_1} \ar[rd]^{p_2}& \\
\underline{\Lambda}^{\alpha}_{\CalI,n,s} & &
 \underline{\Lambda}^{\beta}_{\CalI,n-s,0}
}
\end{equation}
The main theorem of this section is the following, where $\gamma=s[\CalI]$ and $\beta=\alpha-s[\CalI]$:
\begin{thm}\label{CorrespIrr}
We have a natural bijection between irreducible components of
$\underline{\Lambda}^{\alpha}_{\CalI,n,s}$ and irreducible components
of $\underline{\Lambda}^{\beta}_{\CalI,n-s,0}$. If $Z_1 \leftrightarrow Z_2$ under this
correspondence then we have $\dim Z_1 = \dim Z_2 -\langle \beta, \gamma \rangle -\langle \gamma, \beta\rangle-\langle \gamma, \gamma \rangle$. 
\end{thm}
We denote the two applications of the theorem in the following way:
\begin{equation*}
\xymatrix{ \mathrm{Irr}(\underline{\Lambda}^{\alpha}_{\CalI,n,s}) \ar@<1ex>[r]^{f_\CalI^\mathrm{max}} & \mathrm{Irr}(\underline{\Lambda}^{\beta}_{\CalI,n-s,0})\ar@<1ex>[l]^{e_\CalI^s}
}
\end{equation*}

We have to study the maps $p_1$ and $p_2$. We start with $p_2$.
\begin{prop}\label{mapP2}
The map $p_2$ is smooth with connected fibers of dimension
$-\langle \beta, \gamma \rangle -\langle \gamma, \beta\rangle-\langle \gamma, \gamma \rangle$.
\end{prop}
\begin{cor} The map $p_2$ induces a natural bijection between
  irreducible components of
  $\underline{\Lambda}_{\CalI,n-s,0}^{\beta}$ and irreducible components of $\underline{\mathcal{E}}_{\CalI,n,s}^{\alpha}$. Moreover, if $Z_2 \leftrightarrow Z_3$ under
  this correspondence we have $\dim Z_3=\dim Z_2 -\langle \beta, \gamma \rangle -\langle \gamma, \beta\rangle-\langle \gamma, \gamma \rangle$.
\end{cor}
As usual, to prove proposition (\ref{mapP2}) we need to study locally the diagram
(\ref{diag}). To do so, define the locally closed subvariety
$S_{m,(\CalI,n,s)}^{\alpha}$ of $S_m^{\alpha}$ by 
\[
S_{m,(\CalI,n,s)}^{\alpha}=\{ (\phi,\CalF,f) \in S_m^{\alpha}\ | \
(\CalF,f)\in \underline{\Lambda}_{\CalI,n,s}^{\alpha}\}.
\]
For an integer $m<<0$ such that $\CalI$ and $\CalF$ are generated in degree $m$ set $d_1^{(t)}=\langle \mathcal{L}_t(m\vec{c}),\CalI^ s
\rangle(=\dim \Hom(\mathcal{L}_t(m\vec{c})), \CalI^s)$ for $t\in\Pi$, $d'^{(t)}= \langle
\mathcal{L}_t(m\vec{c}),\CalF\rangle(= \dim \Hom(\mathcal{L}_t(m\vec{c}),\CalF))$ and
$d_2^{(t)}=d'^{(t)}-d_1^{(t)}$. Define
\begin{align*}
E^{\alpha,\geq m}_{\CalI,s,n}=\{(\phi,\CalF,f,i,h_1,h_2)| \ (\phi,\CalF,f) \in
S_{m,(\CalI,n,s)}^{\alpha},
i:\CalI^s \hookrightarrow \Ker f ,\ h_1:\Cc^{d_1^{(t)}} \simeq
\Hom(\mathcal{L}_t(m\vec{c}),\CalI^s),\\
h_2:\Cc^{d_2^{(t)}}\simeq \Hom(\mathcal{L}_t(m\vec{c}), \CalF/i(\CalI^s)
\}
\end{align*}
The group $G=\prod_{t\in \Pi} GL_{d'^{(t)}}\times GL_{d_1^{(t)}} \times GL_{d_2^{(t)}}$ acts
naturally on $E^{\alpha,\geq m}_{\CalI,s,n}$ and the quotient stack is
$\underline{\mathcal{E}}_{\CalI,n,s}^{\alpha,\geq m}$.\\
Introduce 
\[C=\{(V^{(t)},a^{(t)},b^{(t)})_{t\in \Pi}\ | \ V^{(t)} \subseteq \Cc^{d'^{(t)}}, \ a^{(t)}:V^{(t)} \simeq \Cc^{d_1^{(t)}}, b^{(t)}: \Cc^{d'^{(t)}}/V^{(t)}\simeq \Cc^{d_2^{(t)}}\}.
\]
We define $q_2$ as follows:
\[
\begin{array}{cccc}
q_2: & E^{\alpha,\geq m}_{\CalI,n,s} & \rightarrow &
 S_{m,(\CalI,n-s,0)}^{\beta}\times S_{m,(\CalI,s,s)}^{\gamma} \times
 C \\
 & (\phi,\CalF,f,i,h_1,h_2) & \mapsto &
 ((\psi_1,\CalG,g),(\psi_2,\CalI^s,0),(V^{(t)},a^{(t)},b^{(t)}))
\end{array}
\]
where:
\begin{enumerate}
\item $\CalG:=\CalF/i(\CalI^s)$,
\item $g:\CalG\rightarrow \CalG(\vec{\omega})$ is induced by $f$,
\item $\psi_1: \bigoplus_{t\in \Pi}\mathcal{L}(m\vec{c})^{d_2^{(t)}} \twoheadrightarrow \CalG$ is deduced from
  $\phi$ and $h_2=(h_2^{(t)})_{t\in\Pi}$,
\item $\psi_2:\bigoplus_{t\in \Pi}\mathcal{L}_t(m\vec{c})^{d_1^{(t)}} \twoheadrightarrow
  \CalI^s$ is deduced from $\phi$ and $h_1=(h_1^{(t)})_{t\in \Pi}$,
\item $(V^{(t)},a^{(t)},b^{(t)})\in Z$ is defined by $V^{(t)}=\phi_{t*}(\Hom
  (\mathcal{L}_t(m\vec{c}),\CalI^s))\subseteq
  \Cc^{d_1^{(t)}+d_2^{(t)}}=\phi_{t*}(\Hom(\mathcal{L}_t(m\vec{c}),\CalF))$ (via $i$) and $a^{(t)}$ and
  $b^{(t)}$ are deduced from the diagram
\begin{equation*}
\xymatrix{
0 \ar[r] & \Hom(\mathcal{L}_t(m\vec{c}),\CalI^s) \ar[r]^i\ar[d]_{h_1^{(t)}}& \Hom(\mathcal{L}_t(m\vec{c}),\CalF)
\ar[r]\ar[d]_{\phi_{t*}}& \Hom(\mathcal{L}_t(m\vec{c}),\CalG) \ar[r]\ar[d]_{h_2^{(t)}} & 0 \\
0 \ar[r]& \Cc^{d_1^{(t)}} \ar[r] &\Cc^{d'^{(t)}} \ar[r] &\Cc^{d_2^{(t)}}\ar[r]& 0
}
\end{equation*}
\end{enumerate}
\begin{lem}
The map $q_2$ is an affine fibration, with fibers of dimension $-\langle
\gamma,\beta\rangle -\langle \beta,\gamma \rangle+ d_1d_2 +s(n-s)$.
\end{lem}
\begin{dem}
Let us fix some notation. We will write $\CalO_\CalI=\bigoplus_{t\in \Pi} \mathcal{L}_t(m\vec{c})^{d_1^{(t)}}$, $\CalO_\CalG=\bigoplus_{t\in \Pi}\mathcal{L}_t(m\vec{c})^{ d_2^{(t)}}$, $\CalO_\CalF=\bigoplus_{t\in\Pi}\mathcal{L}_t(m\vec{c})^{ d'^{(t)}}$.\\
For maps
\begin{equation*}\begin{cases}
\psi_1:\CalO_\CalG\twoheadrightarrow \CalG \\
\psi_2:\CalO_\CalI\twoheadrightarrow \CalI^s \\
\phi:\CalO_\CalF\twoheadrightarrow \CalF 
\end{cases}
\end{equation*}
write $\mathcal{K}_\CalI$, $\mathcal{K}_\CalG$ and
  $\mathcal{K}_\CalF$ for the corresponding kernels. We denote
  $i_\CalI :\mathcal{K}_\CalI \hookrightarrow
  \CalO_\CalI$ and $i_\CalG :\mathcal{K}_\CalG \hookrightarrow
  \CalO_\CalG$ the corresponding injections.\\
Let us describe the morphism $q_2$ in terms of some diagrams. Objects
  in the space $E_{\CalI,n,s}^{\alpha,\geq n}$ are in canonical bijection
  with commutative diagrams
\begin{equation}\label{elementfiber}
\xymatrix{
  0\ar[d] & 0\ar[d] & \\
 \CalO_\CalI \ar[d]^{a'}\ar[r]^{\psi_2}&\CalI^s\ar[r]\ar[d]^i&0 \\
 \CalO_\CalF\ar[d]^{b'}
 \ar[r]^\phi&\mathcal{F}\ar[r]\ar[d] & 0\\
 \CalO_\CalG
 \ar[d]\ar[r]^{\psi_1}&\mathcal{G}\ar[d]\ar[r]&0 \\
 0 & 0 & 
}
\end{equation}
together with a map $f\in \Hom(\CalF,\CalF(\canbun))$ such that
$i:\CalI^s \hookrightarrow \Ker f$. We also require the maps $\psi_1$, $\psi_2$ and $\phi$ respectively induce isomorphisms $\Cc^{d_2^{(t)}}\simeq\Hom(\mathcal{L}_t(m\vec{c}),\CalG)$, $\Cc^{d_1^{(t)}}\simeq\Hom(\mathcal{L}_t(m\vec{c}),\CalI^s)$ and $\Cc^{d'^{(t)}}\simeq\Hom(\mathcal{L}_t(m\vec{c}),\CalF)$.\\
Indeed, the maps $\psi_2,\psi_1$
are deduced from $h_1,h_2$ by the formulas:
\begin{align*}
\psi_2=\mathrm{can} \circ h_1,\\
\psi_1=\mathrm{can} \circ h_2
\end{align*}
where can is the evaluation map and $a',b'$ are defined uniquely in order to make (\ref{elementfiber}) commute. Recall that in the construction of
$\mathrm{Hilb}_{\mathcal{E}_n^\alpha}$, two maps
$\phi:\oplus_{t\in \Pi}\mathcal{L}_t(m)^{d'^{(t)}}\twoheadrightarrow \CalF$, $\phi':\oplus_{t\in\Pi}\mathcal{L}_t(m)^{d'^{(t)}}
\twoheadrightarrow \CalF'$ are equivalent if $\Ker \phi=\Ker
\phi'$. We use the same equivalence relation for diagrams.\\
Similarly, points in $S_{m,(\CalI,n-s,0)}^{\beta}\times S_{m,(\CalI,s,s)}^{\gamma} \times
 C$ correspond bijectively to diagrams
\begin{equation}\label{elementbottom}
\xymatrix{
  0 \ar[d]& & \\
 \CalO_\CalI \ar[d]^{a'}\ar[r]^{\psi_2}&\CalI^s\ar[r]&0 \\
 \CalO_\CalF\ar[d]^{b'}
 & & \\
 \CalO_\CalG
 \ar[d]\ar[r]^{\psi_1}&\mathcal{G}\ar[r]&0 \\
0 & &
}
\end{equation}
together with an element $g\in \Hom(\CalG,\CalG(\canbun))$ (and with the additional condition that $\psi_1$ and $\psi_2$ induce isomorphisms $\Cc^{d_2^{(t)}}\simeq\Hom(\mathcal{L}_t(m\vec{c}),\CalG)$ and $\Cc^{d_1^{(t)}}\simeq\Hom(\mathcal{L}_t(m\vec{c}),\CalI^s)$).\\
The horizontal sequences are the elements of
$S_{m,(\CalI,n-s,0)}^{\beta}$ and $S_{m,(\CalI,s,s)}^{\gamma}$, and the
vertical sequence is deduced from $(V^{(t)},a^{(t)},b^{(t)})$ via
\begin{equation*}
\xymatrix{
0 \ar[r] & V^{(t)} \ar[r] \ar[d]^{a^{(t)}} & \Cc^{d_1^{(t)}+d_2^{(t)}} \ar[r] \ar@2{{}-{}}[d]& \Cc^{d_1^{(t)}+d_2^{(t)}}/V^{(t)} \ar[d]^{b^{(t)}} \ar[r]& 0\\
0 \ar[r] & \Cc^{d_1^{(t)}} \ar[r] & \Cc^{d_1^{(t)}+d_2^{(t)}} \ar[r]& \Cc^{d_2^{(t)}}  \ar[r]& 0
}
\end{equation*}
and then tensoring by
$\mathcal{L}_t(m\vec{c})$.\\
The map $q_2$ assigns to a diagram as in (\ref{elementfiber}) its
subdiagram (\ref{elementbottom}) and the element $g$ deduced from
$f$. We may complete the diagrams (\ref{elementfiber}),
(\ref{elementbottom}) by adding kernels of $\psi_2,\psi_1,\phi$:
\begin{equation}\label{elementfiber2}
\xymatrix{
  & 0\ar[d] & 0\ar[d] & 0\ar[d] & \\
 0\ar[r]& \mathcal{K}_\CalI\ar[r]^{i_\CalI}\ar[d] &
  \CalO_\CalI \ar[d]^{a'}\ar[r]^{\psi_2} &\CalI^s\ar[r]\ar[d]^i&0 \\
 0 \ar[r]& \mathcal{K}_\CalF \ar[r]^{i_\CalF} \ar[d]& \CalO_\CalF\ar[d]^{b'}
 \ar[r]^{\phi}&\mathcal{F}\ar[r]\ar[d] & 0\\
 0\ar[r] & \mathcal{K}_\CalG \ar[r]^{i_\CalG} \ar[d] & \CalO_\CalG
 \ar[d]\ar[r]^{\psi_1}&\mathcal{G}\ar[d]\ar[r]&0 \\
  & 0 &0 & 0 & 
}
\end{equation}
\begin{equation}\label{elementbottom2}
\xymatrix{
  &  & 0\ar[d] &  & \\
 0\ar[r]& \mathcal{K}_\CalI\ar[r]^{i_\CalI} &
  \CalO_\CalI \ar[d]^{a'}\ar[r]^{\psi_2} &\CalI^s\ar[r]&0 \\
  & & \CalO_\CalF\ar[d]^{b'}
 & & \\
  0\ar[r] & \mathcal{K}_\CalG \ar[r]^{i_\CalG}  & \CalO_\CalG
 \ar[d]\ar[r]^{\psi_1}&\mathcal{G}\ar[r]&0 \\
  &  &0 &  & 
}
\end{equation}
Let us fix a point $x=(\psi_1,\CalG,g,\psi_2,\CalI^s,V,a,b)\in
S_{m,(\CalI,n-s,0)}^{\beta}\times S_{m,(\CalI,s,s)}^{\gamma} \times C$, and
denote the fiber $F=q_2^{-1}(x)$.\\
We will use the following lemma:
\begin{lem}
We have a (non-canonical) isomorphism:
\[
\mathcal{K}_\mathcal{F} \simeq \mathcal{K}_\CalG \oplus \mathcal{K}_\CalI
\]
Moreover, any injective morphism $\mathcal{K}_\CalG \oplus \mathcal{K}_\CalI \hookrightarrow \CalO_\CalF$ such that:\\
- the restriction to $\mathcal{K}_\CalI$ is $i_\CalI$,\\
- the induced morphism $\mathcal{K}_\CalG \rightarrow \CalO_\mathcal{G}$ is $i_\CalG$,\\
induces an element $\phi: \CalO_\CalF \rightarrow \CalO_\CalF/\mathcal{K}_\CalG \oplus \mathcal{K}_\CalI$ in $Q_m^\alpha$.
\end{lem}
\begin{dem}
We first apply the functor $\Hom(\mathcal{L}_t(m\vec{c}), \_)$, where $t\in\Pi$, to the short exact sequence:
\[
0\rightarrow \mathcal{K}_\CalI \rightarrow \CalO_\CalI \rightarrow \CalI^s \rightarrow 0
\]
in order to have the long exact sequence:
\begin{equation}
0 \rightarrow \Hom(\mathcal{L}_t(m\vec{c}), \mathcal{K}_\CalI) \rightarrow \Hom(\mathcal{L}_t(m\vec{c}), \mathcal{O}_\CalI) \xrightarrow{\kappa} \Hom(\mathcal{L}_t(m\vec{c}),\CalI^s) \rightarrow \Ext(\mathcal{L}_t(m\vec{c}), \mathcal{K}_\CalI) \rightarrow 0
\end{equation}
But the map $\kappa$ restricted to $\Hom(\mathcal{L}_t(m\vec{c}), \mathcal{L}_t(m\vec{c}))^{d_t(m,\alpha)}\simeq \Cc^{d_t(m,\alpha)}$ is exactly the map $\psi_{2,t*}$, so the condition for $\psi_2$ in $Q_m^\gamma$ implies that $\kappa$ is surjective, hence $\Ext(\mathcal{L}_t(m\vec{c}), \mathcal{K}_\CalI)=0$ for any $t$, which gives that $\Ext(\CalO_\CalG,\mathcal{K}_\CalI)=0$.\\
Now we apply the functor $\Hom(\_,\mathcal{K}_\CalI)$ to the short exact sequence:
\[
0\rightarrow \mathcal{K}_\mathcal{G} \rightarrow \CalO_\mathcal{G} \rightarrow \mathcal{G} \rightarrow 0
\]
which gives a long exact sequence whose end is:
\[
\Ext(\CalO_\CalG,\mathcal{K}_\CalI) \rightarrow \Ext(\mathcal{K}_\CalG,\mathcal{K}_\CalI) \rightarrow 0
\]
which then gives $\Ext(\mathcal{K}_\CalG,\mathcal{K}_\CalI)=0$. We deduce that the short exact sequence
\[
0\rightarrow \mathcal{K}_\CalI \rightarrow \mathcal{K}_\mathcal{F} \rightarrow \mathcal{K}_\mathcal{G} \rightarrow 0
\]
splits, as required.\\
For the second part of the lemma, take a morphism $a:\mathcal{K}_\CalG\oplus \mathcal{K}_\CalI \hookrightarrow \CalO_\CalF$ as in the lemma. Its cokernel $\CalF$ has obviously the right class in $K(\mathrm{Coh}_\mathbb{X})$. We denote by $\phi$ the map $\CalO_\CalF \rightarrow \CalF$ the quotient morphism. It remains to prove that for any $t\in \Pi$, we have:
\[
\phi_{t*}: \Hom(\mathcal{L}_t(m\vec{c}), \CalF) \simeq \Cc^{d^{(t)}}.
\]
But as $\Ext(\mathcal{L}_t(m\vec{c}),\CalI)=0$, we have the following diagram:
\begin{equation*}
\xymatrix{
0 \ar[r] & \Hom(\mathcal{L}_t(m\vec{c}),\CalI) \ar[r] \ar[d] & \Hom(\mathcal{L}_t(m\vec{c}),\CalF) \ar[r] \ar[d]& \Hom(\mathcal{L}_t(m\vec{c}),\CalG) \ar[d]^{b^{(t)}} \ar[r]& 0\\
0 \ar[r] & \Cc^{d_1^{(t)}} \ar[r] & \Cc^{d^{(t)}} \ar[r]& \Cc^{d_2^{(t)}}  \ar[r]& 0
}
\end{equation*}
where the left and right vertical arrows are isomorphisms by hypothesis. The result follows.
\end{dem}
By the preceding lemma, the set of classes of maps $\phi:\CalO_\CalF \twoheadrightarrow \CalF$
making (\ref{elementfiber}) commutative is in bijection with the set of
subsheaves $\mathcal{K}_\CalF \subseteq \CalO_\CalF$ satisfying:
\begin{equation}\label{condkernel}
\begin{cases}
\mathcal{K}_\CalF \cap
a'(\CalO_\CalI)=a'(\mathcal{K}_\CalI)\\
b'(\mathcal{K}_\CalF)=\mathcal{K}_\CalG
\end{cases}
\end{equation}
Subsheaves $\mathcal{K}_\CalF \subseteq \CalO_\CalF$ satisfying
(\ref{condkernel}) form a principal
$\Hom(\mathcal{K}_\CalG,\CalI^s)$-space. Indeed
\[
\{\mathcal{K}_\CalF \subseteq \CalO_\CalF\ |\ \text{(\ref{condkernel}) is
  satisfied}\}= \{ \mathcal{K}_\CalF' \subseteq
  \CalO_\CalF/\mathcal{K}_\CalI\ | \ \mathcal{K}_\CalF'\cap
  \CalI^s=0,b'(\mathcal{K}_\CalF')=\mathcal{K}_\CalG\} 
=\{s: \mathcal{K}_\CalG \rightarrow
  \CalO_\CalF/\mathcal{K}_\CalI\ |\  b'\circ s=Id_{\mathcal{K}_\CalG}\}
\]
and if $s,s'$ are two sections $\mathcal{K}_\CalG \rightarrow
\CalO_\CalF/ \mathcal{K}_\CalI$ as above then the difference $s-s'$ is in $\Hom(\mathcal{K}_\CalG,\CalO_\CalI/\mathcal{K}_\CalI)=\Hom(\mathcal{K}_\CalG,\CalI^s)$.\\
For convenience, let us choose a section $s_0$ as above. This
corresponds to an identification
$\CalO_\CalF/\mathcal{K}_\CalI\simeq \CalI^s \oplus
\CalO_\CalG$. Then to $u\in \Hom(\mathcal{K}_\CalG,\CalI^s)$ we
associate the diagram
\begin{equation}
\xymatrix{
  & 0\ar[d] & 0\ar[d] & 0\ar[d] & \\
 0\ar[r]& \mathcal{K}_\CalI\ar[r]^{i_\CalI}\ar[d] &
  \CalO_\CalI \ar[d]^{a'}\ar[r]^{\psi_2} &\CalI^s\ar[r]\ar[d]^i&0 \\
 0 \ar[r]& \mathcal{K}_\CalF \ar[r]^{i_\CalF} \ar[d]& \CalO_\CalF\ar[d]^{b'}
 \ar[r]^{\phi_u}&\mathcal{F}\ar[r]\ar[d] & 0\\
 0\ar[r] & \mathcal{K}_\CalG \ar[r]^{i_\CalG} \ar[d] & \CalO_\CalG
 \ar[d]\ar[r]^{\psi_1}&\mathcal{G}\ar[d]\ar[r]&0 \\
  & 0 &0 & 0 & 
}
\end{equation}
Note that $\CalF\simeq\mathrm{Coker}(\mathcal{K}_\CalG
\hookrightarrow^{(i_\CalG,u)} \CalO_\CalG \oplus \CalI^s)$.\\
 It remains to describe the possible choices for the map $f$ in the fiber. Such an element
 verifies two conditions:

(*) $f|_{\CalI^s}=0$

(**) $f'|_\CalG=g$, where $f' \in \Hom(\CalG,\CalG(\canbun))$ is deduced
from $f$.\\
We have the short exact sequence derived from $u$:
\[
0 \rightarrow \CalI^s \rightarrow \CalF \rightarrow \CalG
\rightarrow 0
\]
From the following commutative diagram, where vertical maps are
deduced from Serre duality
\begin{equation*}
\xymatrix{
 0 \ar[r]& \Ext(\CalF,\CalG)^* \ar[r] \ar@2{{}-{}}[d] & \Ext(\CalF,\CalF)^* \ar@2{{}-{}}[d]\ar[r]& \Ext(\CalF,\CalI^s)^*\ar@2{{}-{}}[d] \\
0 \ar[r]& \Hom(\CalG,\CalF(\canbun)) \ar[r] & \Hom(\CalF,\CalF(\canbun)) \ar[r]^{|_{\CalI^s}} & \Hom(\CalI^s,\CalF(\canbun))
}
\end{equation*}
We see that $f_{|\CalI^s}=0$ is equivalent to $f\in \Ext(\CalF,\CalG)^*$. The other condition is given by:
\begin{equation*}
\begin{array}{ccccccccc}
0 & \rightarrow & \Ext(\CalI^s,\CalG)^* & \rightarrow & \Ext(\CalF,\CalG)^* & \rightarrow & \Ext(\CalG,\CalG)^* & \xrightarrow{\theta_u} & \Hom(\CalI^s,\CalG)^*\\
 & & & & f & \mapsto & g & &
\end{array}
\end{equation*}
where $\theta_u$ is the connecting morphism.\\
The possible choices of $f$ in the fiber is then a principal $\Ext(\CalI^s,\CalG)^*$-space, when we have the condition $\theta_u(g)=0$.\\
To sum up, we have shown that the fiber $F$ is isomorphic to the
subspace of pairs $(u,v)\in \Hom(\mathcal{K}_\CalG,\CalI^s)\oplus
\Ext(\CalI^s,\CalG)^*$ satisfying $\theta_u(g)=0$. We need to
describe more precisely the map $\theta_u$.
\begin{lem}
The map $\theta_u:\Ext(\CalG,\CalG)^* \rightarrow
\Hom(\CalI^s,\CalG)^*$ is given by
\[
\theta_u(g)(h)= a_g(h\circ u)
\]
for any $h\in \Hom(\CalI^s,\CalG)$, where $a_g$ is the image
of $g$ in $\Hom(\mathcal{K}_\CalG,\CalG)^*$.
\end{lem}
\begin{dem}
We claim that the following diagram is commutative
\begin{equation*}
\xymatrix{
\Ext(\CalG,\CalG)^* \ar[r]^{\theta_u} \ar@{^{(}->}[d] & \Hom(\CalI^s,\CalG)^* \ar@{^{(}->}[d] \\
\Hom(\mathcal{K}_\CalG,\CalG)^* \ar[r]^{\theta_u'} & \Hom(\CalO_\CalG\oplus \CalI^s,\CalG)^*
}
\end{equation*}
where $\theta_u'$ is induced by the injection $\mathcal{K}_\CalG
\xrightarrow{(i_\mathcal{G},u)} \CalO_\CalG \oplus \CalI^s$. \\
To see this, apply $\Hom(.,\CalG)$ to the diagram
\[
\xymatrix{ & & \CalI^s\ar[d]\ar[dr] & & \\
0 \ar[r]&\mathcal{K}_\CalG \ar[r]\ar@{=}[d] &\CalO_\CalG \oplus \CalI^s \ar[r] \ar[d]
& \CalF \ar[r] \ar@{->>}[d] & 0\\
0\ar[r] & \mathcal{K}_\CalG \ar[r] & \CalO_\CalG \ar[r] & \CalG \ar[r]& 0
}
\]
to get the construction of the connecting morphism $\theta^*_u$
\[
\xymatrix{
0 \ar[r] & \Hom(\CalG,\CalG) \ar[r]\ar[d]& \Hom(\CalF,\CalG) \ar[r]\ar[d]
  & \Hom(\CalI^s,\CalG) \ar@{=}[d] \ar@{-->}[dl]\ar `[rr]`[rrddd][dddll]^{\theta_u^*}&  &\\
0 \ar[r] & \Hom(\CalO_\CalG,\CalG) \ar[r]\ar[d] & \Hom(\CalO_\CalG \oplus
\CalI^s,\CalG) \ar@{-->}[dl]\ar@<1ex>[r]^-{proj} \ar[d]^{(i_\mathcal{G},u)} &
\Hom(\CalI^s,\CalG) \ar[r]\ar@<1ex>[l]^-{can} & 0 &\\
 & \Hom(\mathcal{K}_\CalG,\CalG)\ar[d] \ar@{=}[r] &
 \Hom(\mathcal{K}_\CalG,\CalG) & && \\
& \Ext(\CalG,\CalG) & & & &
}
\]
as the composition of the two dotted arrows and the map
$\Hom(\mathcal{K}_\CalG,\CalG) \rightarrow \Ext(\CalG,\CalG)$, which is exactly the dual
of our claim.\\
So for $h\in \Hom(\CalI^s,\CalG)$, we have $\theta_u'(a_g)(h)=a_g(h \circ (i_\CalG,u))$. Then $\theta_u$ is obtained by evaluating $a_g$ on the projection of $h \circ (i_\CalG, u)$ into $\Hom(\CalI^s,\CalG)$, i.e. $\theta_u(g)(h)=a_g(h\circ u)$.
\end{dem}
We can now consider a new linear map $\theta$ defined from $\theta_u$, this time considering the dependence on $u$:
\[
\begin{array}{cccc}
\theta: & \Hom(\mathcal{K}_\CalG,\CalI^s)  & \rightarrow & \Hom(\CalI^s,\CalG)^*\\
 & u & \mapsto & \theta_u(g)
\end{array}
\]
We have proved the following statement:
\begin{lem}\label{fiber}
The fiber $F$ is isomorphic to $\Ext(\CalI^s,\CalG)^* \oplus \Ker \theta$.
\end{lem}
It remains to give the dimension of the fiber $F$. We need the following lemma:
\begin{lem}\label{dim}
We have $(\Im \theta)^\bot= \Hom(\CalI^s,\Ker g)$.
\end{lem}
\begin{dem}
Take $h\in \Hom(\CalI^s,\CalG)$, and define $\CalH:=\Im h$. Now
$h\in (\Im\theta)^\bot$ is equivalent to:
\[
\forall u\in \Hom(\mathcal{K}_\CalG,\CalI^s),\ a_g(h\circ u)=0.
\]
We have a natural map $\Hom(\mathcal{K}_\CalG,\CalI^s) \xrightarrow{p}
\Ext(\CalG,\CalI^s)$, and by definition of $a_g$, we have
$a_g(v)=g(p(v))$ for any $v\in \Hom(\mathcal{K}_\CalG,\CalI^s)$.
\begin{equation*}
\xymatrix{
\Hom(\mathcal{K}_\CalG,\CalI^s) \ar[r]^p \ar[d]^{a_g} &
\Ext(\CalG,\CalI^s) \ar[d]^g \\
\Cc & \Cc
}
\end{equation*}
From the surjection $h: \CalI^s \twoheadrightarrow \CalH$, we
have a surjection $\Ext(\CalG,\CalI^s)\twoheadrightarrow
\Ext(\CalG,\CalH)$, and we can make use of the following commutative
diagram:
\begin{equation*}
\xymatrix{
\Hom(\mathcal{K}_\CalG,\CalI^s) \ar[r]^{h_*'} \ar[d]^{p} &
\Hom(\mathcal{K}_\CalG,\CalH) \ar[d]^{p'} &
\\
\Ext(\CalG,\CalI^s) \ar[d] \ar[r]^{h_*}& \Ext(\CalG,\CalH) \ar[r] & 0\\
 0 & &
}
\end{equation*}
So that we have the following chain of equivalence:
\[
\begin{array}{ccc}
h \in (\Im \theta)^\bot & \Leftrightarrow & \forall u\in
\Hom(\mathcal{K}_\CalG,\CalI^s),\ 0=a_g(h\circ u)=g(p(h \circ u)) \\
& \Leftrightarrow & g(p(h_*'(\Hom(\mathcal{K}_\CalG, \CalI^s))))=0\\
& \Leftrightarrow & g(h_*(p'(\Hom(\mathcal{K}_\CalG,\CalI^s))))=0 \\
& \Leftrightarrow & g(h_*(\Ext(\CalG,\CalI^s)))=0  \mathrm{\  (by\
  surjectivity\ of\ }p')\\
& \Leftrightarrow & g|_{\Ext(\CalG,\CalH)}=0 \mathrm{\  (by\ surjectivity\ of\ }h_*).
\end{array}
\]
Now the restriction $g|_{\Ext(\CalG,\CalH)}$ is equal to the image of
$g$ by the morphism $\Ext(\CalG, \CalG)^* \rightarrow
\Ext(\CalG,\CalH)^*$ which is just the restriction morphism, as we see in the following diagram, where the vertical isomorphisms are given by Serre duality:
\[
\xymatrix{
 \Ext(\CalG,\CalG)^* \ar[r] \ar@2{{}-{}}[d] & \Ext(\CalG,\CalH)^*
 \ar@2{{}-{}}[d]\\
\Hom(\CalG,\CalG(\canbun)) \ar[r]^{|_\CalH} & \Hom(\CalH,\CalG(\canbun))
}
\]
so that $g|_{\Ext(\CalG,\CalH)}=g|_\CalH$, this time considered as an
element of $\Hom(\CalG,\CalG(\canbun))$.\\
We have proved that $h\in(\Im \theta)^\bot \Leftrightarrow g|_\CalH=0$, which
by definition is equivalent to $\CalH \subseteq \Ker g$, and to the statement of lemma \ref{dim}.
\end{dem}
Lemma \ref{fiber} gives that $q_2$ is an affine fibration with
connected fibers (or more precisely, as we work with stacks, a vector bundle stack in the sense of \cite{GPHS}, section 3). Lemma
\ref{dim} allows us to compute the dimension of the fiber. As $\dim
\Hom (\CalI^s,\Ker g)=s \dim \Hom(\CalI,\Ker g)=s(n-s)$, we have:
\[
\begin{array}{ccc}
\dim q_2^{-1}(x)  & = & \dim \Ext(\CalI^s,\CalG)+\dim \Ker \theta\\& 
= &\dim \Ext(\CalI^s,\CalG) + \dim
\Hom(\mathcal{K}_\CalG,\CalI^s)-(\dim \Hom(\CalI^s,\CalG)
-s(n-s))\\
&=&\dim \Ext(\CalI^s,\CalG)- \dim
\Hom(\CalI^s,\CalG)+\dim\Hom(\mathcal{K}_\CalG,\CalI^s) +s(n-s)\\ 
&=&-\langle\CalI^s,\CalG\rangle +\langle
\mathcal{K}_\CalG, \CalI^s\rangle +s(n-s)\\
&=&-\langle \CalI^s,\CalG\rangle +\langle \CalO_\CalG,\CalI^s
\rangle-\langle \CalG,\CalI^s \rangle +s(n-s)\\
&=&-\langle \gamma,\beta\rangle -\langle \beta,\gamma
\rangle +\sum_t d_1^{(t)}d_2^{(t)} +s(n-s).
\end{array}
\]
\end{dem}
As the map $q_2$ is $G$-equivariant, we can pass to the quotient
to obtain a map $q_2'$, which is also an affine fibration with
connected fibers:
\[
q_2':\underline{\mathcal{E}}_{\CalI,n,s}^{\alpha,\geq m} \rightarrow
\underline{\Lambda}_{\CalI,n-s,0}^{\beta,\geq m} \times
\underline{\Lambda}_{\CalI,s,s}^{\gamma,\geq m}\times [C/GL_{d'}]
\]
The variety $C$ is an homogenous $GL_{d'}$-variety, where $GL_{d'}=\prod_t GL_{d'(t)}$, (hence smooth) of
dimension $\sum_t (d'^{(t)})^2-d_1^{(t)}d_2^{(t)}$, so the quotient is a smooth (connected) stack of
dimension $-\sum_t d_1^{(t)}d_2^{(t)}$.\\
We have the following diagram:
\[
\xymatrix{
\underline{\mathcal{E}}_{\CalI,n,s}^{\alpha,\geq m} \ar[rr]^-{q_2'}
\ar[rrd]^{p_2} & & \underline{\Lambda}_{\CalI,n-s,0}^{\beta,\geq m}\times
\underline{\Lambda}_{\CalI,s,s}^{\gamma,\geq m}\times [C/GL_{d'}] \ar[d]^{proj}\\
 &  & \underline{\Lambda}_{\CalI,n-s,0}^{\beta,\geq m}
}
\]
As the stack $\underline{\Lambda}_{\CalI,s,s}^{s[ \CalI ]}$ is smooth
connected of dimension
$-\langle \gamma,\gamma \rangle=-\langle s[\CalI],s[\CalI] \rangle=-s^2$, the morphism $p_2$ is smooth
with connected fibers of dimension $-\langle \beta, \gamma \rangle
-\langle \gamma, \beta \rangle - \langle \gamma,\gamma \rangle+s(n-s)$. The proposition (\ref{mapP2}) follows.

Now we study the map $p_1$.
\begin{prop}\label{mapP1}
There is a natural bijection between between irreducible components $Z_1$ of
$\underline{\Lambda}^{\alpha}_{\CalI,s,n}$ and irreducible components $Z_3$
of $\underline{\mathcal{E}}_{\CalI,n,s}^{\alpha}$. Under this
correspondence we have $\dim Z_3=\dim Z_1 +s(n-s)$.
\end{prop}
\begin{dem}
We enlarge our stack
$\underline{\mathcal{E}}_{\CalI,n,s}^{\alpha}$. Define a stack classifying
isomorphism classes of objects:
\[
\underline{F}_{\CalI,n,s}^{\alpha}=\{(\CalF,f,i) \ |\ 
(\CalF,f)\in\underline{\Lambda}_{\CalI,n,s}^{\alpha} , \ i \in
Gr_s^{\Hom(\CalI,\Ker f)} \}
\]
where $Gr_s^{\Hom(\CalI,\Ker f)}$ is the Grassmannian of $s$-dimensional subspaces of $\Hom(\CalI,\Ker f)$, and where as usual morphisms between objects $(\CalF,f,i)$ and
$(\CalF',f',i')$ are isomorphisms $\psi:\CalF \simeq \CalF'$ such that
the following diagrams commute:
\[
\xymatrix{
 \CalF \ar[r]^\psi \ar[d]^f &  \CalF'\ar[d]^{f'} & \CalI^s
 \ar[r]^i \ar[rd]^{i'} & \Ker f \ar[d]^\psi \\
\CalF(\canbun) \ar[r]^\psi & \CalF'(\canbun) &  & \Ker f'
}
\]
The substack $\underline{\mathcal{E}}_{\CalI,n,s}^{\alpha}$ is easily seen
to be an open dense substack of $\underline{F}_{\CalI,n,s}^{\alpha}$, as
the condition $i$ injective is open
in the irreducible variety $Gr_s^{\Hom(\CalI,\Ker f)}$, and the map
  $p_1$ naturally extends to $\underline{F}_{\CalI,n,s}^{\alpha,\beta}$.\\
Define the stacks
\[\underline{G}_{\CalI,n,s}^{\alpha}=\{ (\CalF,f,i,h)\ | \ (\CalF,f,i)
\in \underline{F}_{\CalI,n,s}^{\alpha}, \ h:\Cc^n\simeq \Hom(\CalI,\Ker
f) \}
\]
and
\[\underline{G'}_{\CalI,n,s}^{\alpha}=\{ (\CalF,f,i,h)| \ i \in
Gr_s^n\}
\]with the rest of the data as in
$\underline{G}_{\CalI,n,s}^{\alpha}$. We have natural maps which lead to the
following commutative diagram:
\[
\xymatrix{
\underline{G}_{\CalI,n,s}^{\alpha} \ar[r]^{\nu_1} \ar[d]_u &
\underline{F}_{\CalI,n,s}^{\alpha} \ar[r]^{p_1} &
\underline{\Lambda}_{\CalI,n,s}^{\alpha} \\
\underline{G'}_{\CalI,n,s}^{\alpha} \ar[r]^{\nu_2} &
\underline{\Lambda}_{\CalI,n,s}^{\alpha}\times Gr_n^s \ar[ur]_{proj} &
}
\]
where we have:
\begin{enumerate}
\item The map $u$ is defined by $u(\CalF,f,i,h)=(\CalF,f,i',h)$ where $i' \in Gr_s^n$ is
  deduced from $i$ via $h$. It is clearly an isomorphism.
\item The maps $\nu_1$ and $\nu_2$, defined by
  $\nu(\CalF,f,i,h)=(\CalF,f,i)$, are $GL_n$ principal bundles.
\end{enumerate}
Consequently, this diagram induces a bijection between
Irr($\underline{F}_{\CalI,s,n}^{\alpha}$) and Irr($\underline{\Lambda}_{\CalI,s,n}^{\alpha}$). But as
$\underline{\mathcal{E}}_{\CalI,s,n}^{\alpha}$ is an open dense substack of
$\underline{F}_{\CalI,s,n}^{\alpha}$, it
also gives a bijection between
Irr($\underline{\mathcal{E}}_{\CalI,s,n}^{\alpha}$) and
Irr($\underline{\Lambda}_{\CalI,s,n}^{\alpha}$). Moreover, under this
correspondence $Z_1 \leftrightarrow Z_2$ we have $\dim Z_1 =\dim Z_2
+s(n-s)$.
\end{dem}
Now Theorem \ref{CorrespIrr} is a consequence of propositions (\ref{mapP2}) and (\ref{mapP1})
\subsection{Definition, consequences and first properties}
The following result is a consequence of Theorem \ref{CorrespIrr}:
\begin{thm}  For any $\alpha \in K^+(\Coh_\mathbb{X})$ and any line bundle $\mathcal{L}$, we have the following:
\begin{itemize}
\item  the space $\underline{\Lambda}_{\mathcal{L},0}^\alpha$ is open in $\underline{\Lambda}_\mathbb{X}^\alpha$,
\item the stack $\underline{\Lambda}_\mathbb{X}^\alpha$ is pure of dimension $-\langle \alpha, \alpha \rangle$.
\end{itemize}
\end{thm}
\begin{dem} We first prove the first point. Let us fix a point $(\CalF,f) \in \underline{\Lambda}_\mathbb{X}^\alpha$. Define a map $\psi_{\CalF,f}$ as follows:
\[
\begin{array}{cccc}
\psi_{\CalF,f}: & \Hom(\mathcal{L},\CalF^\mathrm{vec}) & \rightarrow & \Hom(\mathcal{L},\CalF(\vec{\omega}))\\
 & h & \mapsto & f\circ h.
\end{array}
\]
Now remark that if an element $h\in \Hom(\mathcal{L},\CalF^\mathrm{vec})$ is non-zero, it is an injection, and then $\psi_{\CalF,f}(h)=f|_{\mathcal{L}}$. So the condition that there is no injection $\mathcal{L} \hookrightarrow \Ker f$ is equivalent to the fact that for any injection $h:\mathcal{L} \hookrightarrow \CalF$, the restriction $f|_{\mathcal{L}}=\psi_{\CalF,f}(h)$ is non-zero. We then have the equivalence:
\[
\psi_{\CalF,f} \mathrm{\ is\ injective\ }\ \Longleftrightarrow \ (\CalF,f) \in \underline{\Lambda}_{\mathcal{L},0}^\alpha.
\]
The condition on the left is open, so we have that $\underline{\Lambda}_{\mathcal{L},0}^\alpha$ is open in $\underline{\Lambda}_\mathbb{X}^\alpha$.
\\
Now we prove the second point.
We already know the result for any $\alpha$ such that $\mathrm{rk}(\alpha)=0$ (see Section \ref{IrredTor}). Now we proceed by induction on the rank; take a positive integer $r$ and let us suppose that the theorem is true for any $\beta$ with $\mathrm{rk}(\beta)<r$. Take $\alpha$ with $\mathrm{rk}(\alpha)=r$. We have the following:
\begin{equation}\label{LineBundleDecomposition}
\underline{\Lambda}_\mathbb{X}^{\alpha}= \bigcup_{\mathcal{L} \mathrm{\ line\ bundle},\ s>0}
\underline{\Lambda}_{\mathcal{L},s}^{\alpha}.
\end{equation}
Now consider an irreducible component $Z$ of $\underline{\Lambda}_\mathbb{X}^\alpha$ of dimension $d$. Then there exists a line bundle $\mathcal{L}$ and two positive integers $s$ and $n$ such that $Z \cap \underline{\Lambda}_{\mathcal{L},n,s}^\alpha$ is dense in $Z$. From the diagram (\ref{diag2}), it corresponds to an irreducible component of $\underline{\Lambda}_{\mathcal{L},n,s}^{\alpha-s[\mathcal{L}]}$ of dimension $d'=d+\langle \alpha,s[\mathcal{L}] \rangle +\langle s[\mathcal{L}] ,\alpha\rangle -\langle s[\mathcal{L}] ,s[\mathcal{L}]\rangle$. Its closure $Z'$ in $\underline{\Lambda}_\mathbb{X}^{\alpha-s[\mathcal{L}]}$ is irreducible (but maybe smaller an irreducible component), so by induction hypothesis we have $d'=\dim Z'\leq - \langle \alpha -s[\mathcal{L}],\alpha -s[\mathcal{L}]\rangle$.\\
We want to prove that $d'=-\langle \alpha -s[\mathcal{L}],\alpha -s[\mathcal{L}]\rangle$. For this, consider an irreducible component $\tilde{Z}'$ of $\underline{\Lambda}_{\mathcal{L},0}^{\alpha-s[\mathcal{L}]}$ containing $Z'$. By induction hypothesis and the first part of the theorem, $\tilde{Z}'$ is of dimension $-\langle \alpha -s[\mathcal{L}],\alpha -s[\mathcal{L}]\rangle$. Choose the integer $n'$ such that $\tilde{Z}'\cap \underline{\Lambda}_{\mathcal{L}, n',0}^{\alpha-s[\mathcal{L}]}$ is dense in $\tilde{Z}'$, and use diagram (\ref{diag2}) to obtain a corresponding irreducible component of $\underline{\Lambda}_{\mathcal{L},n'+s,s}^\alpha$. Let us denote by $\tilde{Z}$ the closure of this irreducible component in $\underline{\Lambda}_{\mathcal{L},s}^\alpha$. We have that $\tilde{Z}$ is irreducible of dimension $-\langle \alpha,\alpha\rangle$. We obviously have that $Z\subseteq \tilde{Z}$, and as $Z$ is supposed to be an irreducible component, we have the equality $Z=\tilde{Z}$. Hence we have $d=-\langle \alpha, \alpha \rangle$.
\end{dem}

We now consider the set $\textbf{B}:=\mathrm{Irr}(\underline{\Lambda}_\mathbb{X})$ and endow it with the following structure.\\
First, the decomposition as connected components $\underline{\Lambda}_\mathbb{X}=\bigsqcup_{\alpha \in K^+(\Coh_\mathbb{X})} \underline{\Lambda}_\mathbb{X}^\alpha$ gives rise to a weight map:
\[
wt: \textbf{B} \rightarrow K^+(\Coh_\mathbb{X}).
\]
Given an element $Z\in \textbf{B}$ and an indecomposable rigid element $\CalI$, there is only one stratum $\underline{\Lambda}_{\CalI,n,s}^\alpha$ such that $Z\cap \underline{\Lambda}_{\CalI,n,s}^\alpha$ is dense in $Z$. We define two operators $e_\CalI$ and $f_\CalI$ on the set $\textbf{B}$ the following way:
\[
\begin{array}{ccl}
 e_\CalI(Z)  & =  & e_\CalI^{s+1}(f_\CalI^\mathrm{max}(Z)) \\
 f_\CalI(Z)&   =  &\begin{cases}e_\CalI^{s-1}(f_\CalI^\mathrm{max}(Z)) & \text{ if $s>0$}\\ 0 & \mathrm{\ otherwise}\end{cases}
\end{array}
\]
We have already defined in section \ref{Loop} the function $\epsilon_\CalI: \underline{\Lambda}_\mathbb{X} \rightarrow \Zz^+$, and we define a function, with the same notation, on the set $\textbf{B}$ by taking the generic value on an irreducible component.
Finally define $\phi_\CalI: \textbf{B} \rightarrow \Zz$ by the formula:
\[
\phi_\CalI(Z)=\epsilon_\CalI(Z)+<[\CalI],wt(Z)>
\]
\begin{dfn}
The collection $(\textbf{B}, wt, e_\CalI,f_\CalI,\epsilon_\CalI, \phi_\CalI)$, where $\CalI$ describes the set of indecomposable rigid objects in $\mathrm{Coh}_\mathbb{X}$ is called the \emph{loop crystal} associated to the weighted projective line $\mathbb{X}$.
\end{dfn}
As in the case of usual crystals, we can see this data as a colored graph: the vertices correspond to irreducible components, and we draw an arrow from $Z$ to $Z'$ with color $\mathcal{I}$ if $f_\CalI(Z)=Z'$. The following result is immediate from the construction of the loop crystal.
\begin{prop}
We have the following properties, for any indecomposable rigid object $\CalI$ and any element $Z,Z'$ in $\textbf{B}$:
\begin{enumerate}
\item $wt(e_\CalI(Z))=wt(Z)+[\CalI]$ and $wt(f_\CalI(Z))=wt(Z)-[\CalI]$ if $f_\CalI(Z)\neq 0$,
\item $\epsilon_\CalI(e_\CalI(Z))=\epsilon_\CalI(Z)+1$ and $\phi_\CalI(e_\CalI(Z))=\phi_\CalI(Z)-1$,
\item $\epsilon_\CalI(f_\CalI(Z))=\epsilon_\CalI(Z)-1$  (resp. $\phi_\CalI(f_\CalI(Z))=\phi_\CalI(Z)+1$) if $\epsilon_\CalI(Z)\neq 0$,
\item $Z'=f_\CalI(Z)$ if and only if $e_\CalI(Z')=Z$.
\end{enumerate}
\end{prop}
We can also prove the following:
\begin{prop}
The loop crystal graph is connected.
\end{prop}
\begin{dem}
First remark that if $Z$ is an irreducible component of class $\alpha$ with $\alpha$ of rank strictly positive, then there exist a line bundle $\mathcal{L}$ such that $f_\mathcal{L}(Z)$ is non zero by (\ref{LineBundleDecomposition}). By an easy induction on the rank, any irreducible component $Z$ is connected to an irreducible component $Z'$ of rank $0$.\\
Now from the description of the irreducible components of rank $0$, we can apply the operators $f_\CalI$ for some exceptional torsion sheaves $\CalI$ to reduce ourselves to irreducible components of type $Z_\lambda$, for a partition $\lambda$ (as defined in section \ref{Ordinary}). We here use the well-known fact that the crystal of each cyclic quiver is connected, and observe that the $f_\CalI$ are the usual crystal operators via the correspondence seen in section \ref{Except}. Introduce the conjugate partition $\mu$ such that $\mu_i=|\{j\ |\ \lambda_j\geq i\}|$. Let $k=\mathrm{length}(\mu)$. Introduce the following sequence of partitions:
\[
\mu^{(i)}=(\mu_1,\mu_2,\cdots, \mu_i)\ , \ 1\leq i \leq k,
\]
and $\lambda^{(i)}$ the conjugate of $\mu^{(i)}$. For $l\Nn$, define the following vector bundle on $\mathbb{X}$:
\[
\mathcal{V}_l:= \oplus_{k=1}^{l} \CalO(2k\vec{c}).
\]
Consider the space:
\[
Z_l=\{(\mathcal{V}_l,f)\ | \ f\in \Hom(\mathcal{V},\mathcal{V}(\vec{\omega}))\mathrm{\ such\ that\ } \CalO(2k\vec{c})\rightarrow \CalO(2(k+1)\vec{c}+\vec{\omega})\mathrm{\ is\ non-zero\ for\ }1\leq k \leq l-1\}.
\]
Then the closure $\overline{Z}_l$ is an irreducible component of $\underline{\Lambda}_\mathbb{X}$ as it is the closure of the conormal to the substack parametrizing sheaves isomorphic to $\mathcal{V}_l$ (in the tangent stack $T^*\underline{\mathrm{Coh}}_\mathbb{X}$). Indeed the nilpotency condition is automatically checked for any element $f\in \Hom (\mathcal{V}_l,\mathcal{V}_l(\vec{\omega}))$: as there is no non-zero map from $\CalO(2(k+1)\vec{c})$ to any $\CalO(2(j+1)\vec{c}+\omega)$ for $j<k$, the kernel should contain $\CalO(2l\vec{c})$, and by an easy induction $f^l=0$. Now we introduce for a partition $\nu$ and a positive integer $l$ the following space:
\[
Z_{l,\nu}:= \{(\mathcal{V}\oplus \tau_\nu,f)\ | \ (\mathcal{V},f|_{\mathcal{V}})\in Z_l,\ \tau_\nu =\bigoplus_i \CalO_{y_i}^{(\nu_i)},\ \Im f \cap \tau_\nu(\vec{\omega})=\tau_\nu(\vec{\omega})\},
\]
where $y_i$ are distinct ordinary points of $\mathbb{X}$.\\
Via the description of irreducible components of $\underline{\Lambda}_\mathbb{X}$ in \ref{buntor}, the closure $\overline{Z}_{l,\nu}$ corresponds to a conormal bundle to a substack parametrizing vector bundles isomorphic to $\mathcal{V}_l$ in the vector bundle part and the partition $\nu$ for the ordinary torsion sheaf part. It is thus an irreducible component of $\underline{\Lambda}_\mathbb{X}$.
\begin{lem}\label{connected lemma}
For any $1\leq j \leq k$ we have the following:
\[
f_{\CalO((2k-\sum_{m=1}^j\lambda_m)\vec{c})}(\overline{Z}_{k-j+1,\lambda^{(j-1)}})=\overline{Z}_{k-j,\lambda^{(j)}}
\]
and
\[
f_{\CalO((2k-j)\vec{c})}(\overline{Z}_{k-j+1})=\overline{Z}_{k-j}
\]
\end{lem}
\begin{dem}
A generic element of the irreducible component $\overline{Z}_{k-j+1,\lambda^{(j-1)}}$ is of the form $(\mathcal{V}_{k-j+1}\oplus \tau_{\lambda^{(j-1)}}, f)$, where:
\begin{itemize}
\item the kernel of $f$ is of rank $1$, and the restriction of $f$ to $\CalO(m\vec{c})$, $m<k-j+1$, is injective,
\item $\Im f \cap \tau_{\lambda^{(j-1)}}(\vec{\omega}) =\tau_{\lambda^{(j-1)}}(\vec{\omega})$.
\end{itemize}
It follows that the image of $f$ is exactly $(\mathcal{V}_{k-j}\oplus \tau_{\lambda^{(j-1)}})(\vec{\omega})$. Generically, the torsion part of the kernel of $f$ is the socle of $\tau_{\lambda^{(j-1)}}$, i.e. $\bigoplus_i \CalO_{y_i}^{(1)}$. We then deduce that generically the kernel of $f$ is:
\[
\ker f=\CalO((2(k-j+1)-l(\lambda^{(j-1)}))\vec{c}) \oplus \oplus_i\CalO_{y_i}.
\]
Then in the case $l(\lambda)=0$ we have obviously $f_{\CalO((2k-j)\vec{c})}(\overline{Z}_{k-j+1})=\overline{Z}_{k-j}$, as the only injection $\CalO((2k-j)\vec{c}) \hookrightarrow \Ker f$ is an isomorphism into the locally free part of $\Ker f$. Then $\mathcal{V}_{k-j+1}/\CalO(2(k-j))=\mathcal{V}_{k-j}$.\\
In the case $l(\lambda) > 0$, first compute the quotient $\CalF/\Ker f^\mathrm{vec}=\CalF/\CalO((2(k-j+1)-l(\lambda^{(j-1)}))\vec{c})$. For simplicity, define $\sigma=\lambda^{(j-1)}$. We have the exact sequence 
\[
0\rightarrow \Ker f \rightarrow \CalF \rightarrow \CalF/\Ker f \rightarrow 0.
\]
By taking the quotient by $\Ker f^\mathrm{vec}$, we have:
\[
0 \rightarrow \Ker f/\Ker f^\mathrm{vec} \rightarrow \CalF/\Ker f^\mathrm{vec} \rightarrow \CalF/\Ker f \rightarrow 0.
\]
But $\Ker f/ \Ker f^\mathrm{vec}$ is exactly $\Ker f^\mathrm{tor}= \mathrm{soc}(\tau_{\sigma})$. If we decompose $\CalF/\Ker f^\mathrm{vec}$ into its locally free part and its torsion part as $\mathcal{V}_{k-j}\oplus \tau$, we have by restricting to torsion parts:
\[
0 \rightarrow \mathrm{soc}(\tau_{\sigma}) \rightarrow \tau \rightarrow \tau_{\sigma} \rightarrow 0.
\]
We deduce that $\tau$ has support $\{ y_i\}$, with multiplicities $\sigma_i+1$. But as $Z_{l,\sigma}$ is a dense subset of an irreducible component, the objects $(\CalF/\CalO((2(k-j+1)-l(\sigma))\vec{c}), \tilde{f})$ (where $\tilde{f}$ is deduced from $f$) describe, via the diagram (\ref{diag2}), an open dense subset of $f_{\CalO((2(k-j+1)-l(\sigma))\vec{c})}(\overline{Z}_{l,\sigma})$. By the description of the torsion part of irreducible component given in Section (\ref{irredpart}), this implies that generically $\tau=\bigoplus_i \CalO_{y_i}^{(\sigma_i+1)}$.\\
Now take an injection $I :\CalO(a\vec{c})\hookrightarrow \Ker f$. It factors through $\Ker f^\mathrm{vec}$, and generically the quotient $\Ker f^\mathrm{vec}/\CalO(a\vec{c})$ is of the form $\CalO_{z_i}^{(1)}$, where $z_i$ are distinct ordinary points, different from the $y_i$'s. This implies that generically the quotient $\CalF/\CalO(a\vec{c})$ is of the form
\[
\mathcal{V}_{k-j}\oplus  \bigoplus_{i=1}^{|\sigma|} \CalO_{y_i}^{(\sigma_i+1)} \oplus \bigoplus_{i=1}^{|\lambda^{(j)}|-|\lambda^{(j-1)}|} \CalO_{z_i}^{(1)}.
\]
The space of such objects form a dense subspace of the irreducible component $(\mathcal{V}_{k-j},\lambda^{(j)})$, which proves the second part of the lemma.
\end{dem}
In particular, as we have $\overline{Z}_{k,\lambda^{(0)}}=\overline{Z}_k$, $\overline{Z}_{0,\lambda^{(k)}}=Z_\lambda$ and $\overline{Z}_0=\emptyset$, the irreducible components $\emptyset$ and $Z_\lambda$ are connected.
\end{dem}

\subsection{Examples}
In this subsection we give some exemples of computations and describe the irreducible components of the locally free part in the cases $g_\mathbb{X}<1$ and $g_\mathbb{X}=1$.\\

\textbf{Case of the projective line $\mathbb{X}=\P1$}. The corresponding root system is the one of $\widehat{sl}_2$. The irreducible components are indexed by pairs $(\mathcal{V}, \lambda)$, where $\mathcal{V}$ is a vector bundle and $\lambda$ a partition. The vector bundles split into sums of line bundles, and the only indecomposable rigid coherent sheaves are the line bundles $\CalO(k)$, $k\in \Zz$.\\
We have the obvious computation:
\[
f_\CalO^\mathrm{max}((\CalO^n,0))=f_\CalO^n((\CalO^n,0))=\emptyset
\]
from which we deduce that:
\[
f_\CalO((\CalO^n,0))=(\CalO^{n-1},0)\mathrm{\ and\ } e_\CalO((\CalO^n,0))=(\CalO^{n+1},0).
\]
Now a generic quotient $\CalO^n/\CalO(-1)^n$ is of the form $\bigoplus_{i=1}^n \CalO_{y_i}^{(1)}$, where $y_i$ are distinct points of $\P1$, i.e. corresponding to the partition $(1^n)$. Then we have:
\[
f_{\CalO(-1)}^\mathrm{max}((\CalO^n,0))=f_{\CalO(-1)}^n((\CalO^n,0))=(0,(1^n)),
\]
and more generally for $n,l \geq 0$
\[
f_{\CalO(-1)}^\mathrm{max}((\CalO(1)^l\oplus\CalO^n,0))=f_{\CalO(-1)}^{n+l}((\CalO(1)^l\oplus\CalO^n,0))=(0,(1^{n+2l}))
\]
This gives us
\[
f_{\CalO(-1)}((\CalO(1)^l\oplus\CalO^n,0))=(\CalO(1)^{l+1}\oplus \CalO^{n-2},0) \quad \mathrm{for\ }n\geq 2,
\]
\[
f_{\CalO(-1)}((\CalO(1)^l\oplus \CalO,0))=(\CalO(2) \oplus \CalO(1)^{l-1},0)\quad \mathrm{for\ }l\geq 1,
\]
and
\[
f_{\CalO(-1)}((\CalO,0))=(0,(1)).
\]
We can sum up these in the following, where we only draw arrows corresponding to operators $f_\CalO$ and $e_{\CalO(1)}$, and some vertices corresponding to \textit{stable} vector bundles, i.e. of the form $\CalO(k)^a\oplus \CalO(k-1)^b$, for some integer $k$ and some non-negative integers $a$ and $b$.
\[
\xymatrix{
&  \cdots & \cdots & \\
\cdots \ar[r]^{f_\CalO} &\CalO(1)^2 \ar[r]^{f_\CalO} & \CalO(2) \ar[r]^{f_\CalO} \ar[lu]^{e_{\CalO(1)}}& (1^2)\ar[lu]^{e_{\CalO(1)}}\\
\cdots \ar[r]^{f_\CalO} &\CalO(1)\oplus \CalO \ar[r]^{f_\CalO} & \CalO(1) \ar[r]^{f_\CalO} \ar[lu]^{e_{\CalO(1)}}& (1)\ar[lu]^{e_{\CalO(1)}}\\
\cdots \ar[r]^{f_\CalO} &\CalO^2 \ar[r]^{f_\CalO} & \CalO \ar[r]^{f_\CalO} \ar[lu]^{e_{\CalO(1)}} & \emptyset \ar[lu]^{e_{\CalO(1)}}
}
\]
Here for simplicity, for a locally free irreducible component $(\mathcal{V},0)$ we wrote just the vector bundle $\mathcal{V}$, and for a torsion component $(0,\tau)$ we wrote $(\tau)$ .\\
As for non-stable vector bundles, we also proved during lemma \ref{connected lemma} that we have for instance:
\[
\xymatrix{
\CalO(4) \oplus \CalO(2) \oplus \CalO \ar[r]^{f_{\CalO(4)}} \ar[dr]^{f_{\CalO(3)}}  & \CalO(2) \oplus \CalO \ar[r]^{f_{\CalO(2)}} & \CalO\ar[r]^{f_\CalO}  \ar[dr]^{f_{\CalO(-1)}} \ar[ddr]_{f_{\CalO(-2)}} & \emptyset \\
 & (\CalO(2) \oplus \CalO,(1))\ar[dr]^{f_{\CalO(1)}}\ar[ddr]_{f_{\CalO}} & & (1) \\
 & & (\CalO,(2)) \ar[dr]^{f_{\CalO(-1)}}& (1^2)\\
 & & (\CalO,(2,1)) \ar[dr]_{f_{\CalO(-2)}}& (3)\\
 & & & (3,2)
}
\]
\textbf{Case $g_\mathbb{X}<1$}. In this case we have $p\vec{\omega} < 0$, and it easily implies that any $f\in \Hom(\mathcal{V},\mathcal{V}(\vec{\omega}))$, where $\mathcal{V}$ is a vector bundle, is nilpotent. The closure of any conormal bundle to the substack classifying sheaves isomorphic to a given vector bundle lies inside the global nilpotent cone, hence the irreducible components of $\underline{\Lambda}_\mathbb{X}$ are indexed by pairs $(\mathcal{V},\tau)$, where $\mathcal{V}$ is a vector bundle and $\tau$ is an irreducible componentof the torsion part.\\
The Kac-Moody algebra corresponding to the star-shaped diagram (\ref{quiver}) is of finite Dynkin type. Its loop-Kac-Moody algebra is then an affine Lie algebra, corresponding to an affine Dynkin type.\\
The indecomposable rigid objects are of two kinds:
\begin{enumerate}
\item indecomposable vector bundles $\mathcal{V}$. Following Crawley-Boevey (\cite{CB2}), there is a indecomposable vector bundle of class $\alpha$ if and only if $\alpha$ is a positive root, and it is unique up to isomorphism. For instance in type $A$ (i.e. when $n\leq 2$), they are all line bundles, and as soon as $n>2$ there are indecomposable vector bundles of rank $>1$.
\item indecomposable torsion sheaves. These are parametrized by positive roots of rank $0$. The rigid ones are parametrized by \textit{real} roots.
\end{enumerate}
\textbf{Case $g_\mathbb{X}=1$}. In this case the nilpotency condition is non-empty, as $p\vec{\omega}=0$. The corresponding Kac-Moody algebra is of affine Dynkin type, so the loop Kac-Moody algebra is a double-affine, or elliptic, Lie algebra. Unlike the previous case, some positive roots corresponding to (indecomposable) vector bundles are imaginary, meaning that these vector bundles are not rigid. Nevertheless, it is possible to describe explicitly the irreducible components. For this we need to introduce more notation. We will follow the one used in \cite{Sc3}.\\
First define the group homomorphism  $\partial : L(\barp) \rightarrow \Zz$ by $\partial (\vec{x}_i)=\frac{p}{p_i}$. Then we define the \textit{degree map} $\textbf{d}$ as the map of groups
\[
\textbf{d}:K(\mathrm{Coh}_\mathbb{X}) \rightarrow \Zz
\]
defined on generators by $\textbf{d}(\CalO_\mathbb{X}(\vec{x}))=\partial (\vec{x})$.\\
Define the \textit{slope} $\nu(\CalF)$ of a coherent sheaf $\CalF$ as follows
\[
\begin{array}{cccc}
\nu: 	& \mathrm{Coh}_\mathbb{X} & \rightarrow & \Zz \cup \{\infty\}\\
 & \CalF & \mapsto & \textbf{d}(\CalF)/\mathrm{rk}(\CalF).
\end{array}
\]
Here $\nu(\CalF)$ is defined to be $\infty$ when $\CalF$ is a torsion sheaf.\\
We need to introduce the \textit{Harder-Narasimhan filtration} of a coherent sheaf, defined in \cite{GL}.
We say that a coherent sheaf $\CalF$ is \textit{semi-stable} if for any $\CalG \subseteq \CalF$ we have $\nu(\CalG) \leq \nu(\CalG)$. We also define for $q\in \Qq$ the category $\mathcal{C}_q$ to be the full subcategory of $\mathrm{Coh}_\mathbb{X}$ consisiting of the zero sheaf together with sheaves of slope $q$. The following properties are proved in \cite{GL} (for the first three ones) and \cite{LM1} (for the last one):\\
(1) for any $q \in \Qq$, the category $\mathcal{C}_q$ is abelian and closed under extensions,\\
(2) if $\CalF \in \mathcal{C}_q$, $\CalG \in \mathcal{C}_{q'}$ and $q>q'$ then $\Hom(\CalF,\CalG)=0$,\\
(3) (using $g_\mathbb{X}=1$) if $\CalF \in \mathcal{C}_q$, $\CalG \in \mathcal{C}_{q'}$ and $q<q'$ then $\Ext(\CalF,\CalG)=0$,\\
(4) for any $q\in \Qq$ the category $\mathcal{C}_q$ is naturally isomorphic to $\mathcal{C}_\infty$.\\
Any coherent sheaf $\CalF$ has a unique filtration, the Harder-Narasimhan filtration, 
\[
0 \subseteq \CalF_1 \subseteq \cdots \subseteq \CalF_l=\CalF,
\]
such that $\CalF_i/\CalF_{i-1}$ is semistable of slope $\nu_i$ and $\nu_1>\nu_2> \cdots >\nu_l$. The sequence of elements in $K^+(\mathbb{X})$:
\[
\mathrm{HN}(\CalF)= ([\CalF/\CalF_{l-1}],\cdots,[\CalF_2/\CalF_1],[\CalF_1])
\]
is called the \textit{HN-type} of $\CalF$. For any sequence $\underline{\alpha}=(\alpha_1,\alpha_2, \cdots, \alpha_k)$ of elements of $K^+(\mathbb{X})$, we can then define
\[
\mathrm{HN}^{-1}(\underline{\alpha})=\{ \CalF \in \Coh_\mathbb{X}\ | \ \mathrm{HN}(\CalF)=\underline{\alpha}\},
\]
a constructible substack  of $\Coh_\mathbb{X}$. Set $|\underline{\alpha}| :=\sum_i \alpha_i$. We define
\[
T_{\mathrm{HN}^{-1}(\underline{\alpha})}^*:=\{ (\CalF,f) \in T^*\Coh_\mathbb{X}\ | \ \mathrm{HN}(\CalF)=\underline{\alpha}\}
\]
a constructible substack of $T\Coh_\mathbb{X}^*$, and 
\[
\underline{\Lambda}_{\mathrm{HN}^{-1}(\underline{\alpha})}:=\underline{\Lambda}_\mathbb{X} \cap T_{\mathrm{HN}^{-1}(\underline{\alpha})}^*
\]
a constructible substack of $\underline{\Lambda}_\mathbb{X}$. We also define, for any $\alpha\in K^+(\mathbb{X})$
\[
\underline{\Lambda}_\mathbb{X}^{(\alpha)}= \{ (\CalF,f) \in \underline{\Lambda}_\mathbb{X}^\alpha\ |\ \CalF \mathrm{\ semistable\ of\ slope\ } \nu(\alpha) \}.
\]
The proof of the following proposition is similar to the one of lemma \ref{buntor}.
\begin{prop}
The natural morphism
\[
P_{\underline{\alpha}}: \underline{\Lambda}_{\mathrm{HN}^{-1}(\underline{\alpha})} \rightarrow \prod_{i=1}^k \underline{\Lambda}_\mathbb{X}^{(\alpha_i)}
\]
given by $P_{\underline{\alpha}}(\CalF,f)=((\CalF/\CalF_{l-1},f|_{\CalF/\CalF_{l-1}}), \cdots, (\CalF_2/\CalF_{1},f|_{\CalF_2/\CalF_{1}}), (\CalF_1,f|_{\CalF_1}))$ is an affine fibration, hence it induces a correspondence between irreducible components.
\end{prop}
We can now describe all irreducible components of $\underline{\Lambda}_\mathbb{X}^\alpha$: from the previous proposition we have a bijection
\[
\mathrm{Irr}(\underline{\Lambda}_{\mathrm{HN}^{-1}(\underline{\alpha})}) \leftrightarrow \prod_{i=1}^k\mathrm{Irr} (\underline{\Lambda}_\mathbb{X}^{(\alpha_i)})
\]
As the substacks $\underline{\Lambda}_{\mathrm{HN}^{-1}(\underline{\alpha})}$ stratify $\underline{\Lambda}_\mathbb{X}^\alpha$, we have further bijections:
\[
\mathrm{Irr}(\underline{\Lambda}_\mathbb{X}^\alpha) \leftrightarrow \bigcup_{|\underline{\alpha}|=\alpha} \mathrm{Irr}(\underline{\Lambda}_{\mathrm{HN}^{-1}(\underline{\alpha})}) \leftrightarrow \bigcup_{\sum_{i=1}^l\alpha_i=\alpha} \prod_{i=1}^k\mathrm{Irr} (\underline{\Lambda}_\mathbb{X}^{(\alpha_i)})
\]
We can use Property (4) in the previous discussion to determine the irreducible component of $\underline{\Lambda}_\mathbb{X}^{(\alpha_i)}$. The group of automorphism $\Aut(K^+(\mathbb{X}))$ has been determined in \cite{LM2}: it is an extension of $PSL_2(\Cc)$ with the product of the Picard group $\mathrm{Pic}_0(\mathbb{X})$ and a finite group. For any $\alpha\in K^+(\mathbb{X})$ there is an element $\gamma \in \Aut(K^+(\mathbb{X}))$ such that $\nu(\gamma (\alpha))=\infty$; as proved in \cite{LM2}, this action lifts to a natural isomorphism from $\mathcal{C}_{\nu(\alpha)}^\alpha$, the stack classifying semi-stable sheaves of class $\alpha$, to $\mathcal{C}_\infty^{\gamma(\alpha)}$. We thus have an isomorphism:
\[
\underline{\Lambda}_\mathbb{X}^{(\alpha)} \simeq \underline{\Lambda}_\mathbb{X}^{\gamma(\alpha)}
\]
and the irreducible components of $\underline{\Lambda}_\mathbb{X}^{\gamma(\alpha)}$ are known from Section \ref{irredpart}.
\subsection{Principal subspace and representation theory}
In \cite{Sc2}, Schiffmann constructs a canonical basis $\widehat{\bf{B}}$ of the Hall algebra of the category $\text{Coh}_\mathbb{X}$, and proves in the case where the genus $g_\mathbb{X}$ is less or equal to $1$ that the Hall algebra is isomorphic to the quantum version of the algebra $\widehat{U}(\mathcal{L}\mathfrak{n})$. When the algebra $U(\mathcal{L}\mathfrak{g})$ is a Kac-Moody algebra, i.e. in the finite type case, it is conjectured that these bases $\widehat{\bf{B}}$ are compatible with integrable highest weight representations, just like the usual canonical bases $\bf{B}$, despite the fact that they do not correspond to the same positive parts of the enveloping algebras. In the simplest case, i.e. $\mathbb{X}=\mathbb{P}^1$, it is shown that the canonical basis is compatible with the principal subspace $W_0$ of the integrable highest weight $L_{\Lambda_0}$. Here $\Lambda_0$ is the weight of $\widehat{sl}_2$ defined by $\langle \Lambda_0,h_0\rangle =1$ and $\langle \Lambda_0,h_1 \rangle=0$. If we define $\mathfrak{n}$ to be the subalgebra of $\widehat{sl}_2$ generated by $e\otimes t^k$ for $k\in \Zz$, then the principal subspace is
\[
W_0=U(\mathfrak{n}).v_{\Lambda_0},
\]
where $v_{\Lambda_0}$ is the highest weight vector of $L_{\Lambda_0}$. It is also known that the entire space $L_{\Lambda_0}$ is obtained as limit of "twisted" principal subspaces $W_0\subset W_{-1} \subset \cdots $. These twisted principal subspaces are just the images of $W_0$ under the action of the translation lattice $\mathbb{Z}$ in the affine braid group of type $A_1$ on the highest weight module $L_{\Lambda_0}$. The basis $\widehat{\bf{B}}$ is then compatible with $L_{\Lambda_0}$, and is conjectured to be the same as the canonical basis $\bf{B}$.\\
It is also conjectured in \cite{Sc2} that the basis $\widehat{\bf{B}}$ is compatible with the Fourier coefficients of the infinite series, for $l\in \Nn$:
\[
:E(z)^l:=\sum_t \sum_{\substack{l_1t_1+l_2t_2+\cdots l_rt_r\\t_1<\cdots <t_r,l_1+\cdots+l_r=t}}:E_{t_1}^{l_1}\cdots E_{t_r}^{l_r} :,
\]
where $:E_{t_1}^{l_1}\cdots E_{t_r}^{l_r} :=v^{x(l,t)}E_{t_1}^{(l_1)}\cdots E_{t_r}^{(l_r)}$, and $x(l,t)=l^2-\sum_{i\leq j} l_il_j(t_j-t_i+1)$.\\
This suggests that the loop crystal structure obtained should be related to these "semi-infinite" descriptions of integrable highest weight modules, as depicted in \cite{FS}.

\end{document}